\documentclass[12pt]{article}
\usepackage{graphicx}
\usepackage{amsmath}
\usepackage{graphicx}
\usepackage{amsfonts}
\usepackage{mathrsfs}
\usepackage{amssymb}
\usepackage{cite}
\usepackage[multiple]{footmisc}
\usepackage{hyperref}
\usepackage{amsmath, amssymb ,amsthm, amsfonts, amsgen}
\usepackage{xcolor}
\DeclareGraphicsExtensions{.eps,.bmp,.jpg,.pdf,.mps,.png,.gif}
\numberwithin{equation}{section} 
\newcommand{\CC}{\mathbb{C}}
\newcommand{\NN}{\mathbb{N}}

\newcommand{\RR}{\mathbb{R}}
\newcommand{\ZZ}{\mathbb{Z}}

\newtheorem{definition}{\sc Definition}[section]
\newtheorem{teo}{\sc Theorem}[section]
\newtheorem{prop}{\sc Proposition}[section]
\newtheorem{lemma}{\sc Lemma }[section]

\newtheorem{obs}{\sc Remark}[section]

\usepackage[titletoc]{appendix}
\newcommand{\dem}{\noindent {\it Proof.} \mbox{}}    
\newcommand{\lqqd}{\hfill $\square$}      

\begin{document}
\title{Bending shell theories for multiscale materials from $3D$ nonlinear elasticity}

\author{Tiziana Durante\footnote{
Dipartimento di Ingegneria dell'Informazione ed Elettrica e Matematica Applicata, Universit\`a degli Studi di Salerno,Via Giovanni Paolo II, 132 Salerno, Italy. E-mail: tdurante@unisa.it}
\and Luisa Faella\footnote{Dipartimento di Ingegneria Elettrica e dell'Informazione ``M. Scarano",
Universit\`a degli Studi di Cassino e del Lazio Meridionale, Via G. Di Biasio n.43, 03043 Cassino (FR), Italy. E-mail: l.faella@unicas.it} \footnote{European University of Technology EUt+, European Union} \and Pedro Hern\'{a}ndez-Llanos\footnote{
 Instituto de Ciencias de la Ingenier\'ia, Universidad de O'Higgins, Avenida Libertador Bernardo O'Higgins 611, 2841935 Rancagua, Chile. E-mail: pedro.hernandez@postdoc.uoh.cl {\bf(corresponding author)}} \and Ravi Prakash\footnote{Departamento de Matem\'aticas, Facultad de Ciencias F\'isicas y Matem\'aticas, Universidad de Concepci\'on, Avenida Esteban Iturra s/n, Barrio Universitario, Casilla 160-C, Concepci\'on, Chile. Email: rprakash@udec.cl}}
\date{\today}

\maketitle{}

\begin{abstract}

This article derives homogenized bending shell theories starting from three-dimensional nonlinear elasticity.
 The original three-dimensional model contains three small parameters: the two homogenization scales $\varepsilon$ and $\varepsilon^2$ of the material properties and the thickness $h$ of the shell. We obtain different limiting behaviors depending on the limit of various ratios of these three parameters.


\end{abstract}

\section{Introduction}

 Since the early 1990s, the search for lower-dimensional models describing thin three-dimensional structures has attracted considerable interest due to their relevance in the design and manufacturing of new materials with tailored properties. Nonlinear models for plates and shells made of homogeneous materials have been rigorously derived using $\Gamma$-convergence techniques \cite{REF6}, starting from three-dimensional nonlinear elasticity. A hierarchy of limit models has been established through $\Gamma$-convergence, depending on the scaling of the elastic energy $J^{h}$ with respect to the thickness parameter $h \in (0,1]$. Specifically, the scalings $J^h \sim 1$, $J^h \sim h^2$, and $J^h \sim h^4$ correspond to the membrane, bending, and von Kármán regimes, respectively. The first results in this direction for the membrane regime were obtained in \cite{Ledretraoult1} and \cite{Ledretraoult2} for plates and shells, respectively. The nonlinear bending theory for plates was derived in \cite{FJM02}, while the corresponding shell theory was studied in \cite{FJM03}. For the nonlinear von Kármán theory of plates, see \cite{FJM06}. A generalization of the nonlinear von Kármán theory for shells, covering more general scaling laws such as $J^h \sim h^\beta$ with $\beta \geq 4$, can be found in \cite{Lewicakamorapa}.Recently, models of homogenized bending plates and shells for heterogeneous materials have been derived through the simultaneous application of homogenization and dimension reduction techniques, particularly in the case where the relationship between the body’s thickness \( h \) and the material’s oscillations \( \varepsilon(h) \) is characterized by the existence of the limit
\begin{equation}
\gamma := \lim_{h \to 0} \frac{h}{\varepsilon(h)} \in [0, +\infty],
\end{equation}
under suitable assumptions on the stored energy density. Different limit models arise depending on the value of \( \gamma \) (\( \gamma = 0 \), \( \gamma \in (0, +\infty) \), and \( \gamma = +\infty \)); see, for example, \cite{HORVEL14,HORVEL15,HORVEL18} and \cite{Neukammvelcic}. In \cite{BUFDAVFON15}, the authors consider the energy scaling corresponding to the bending plate theory, introducing two distinct material oscillation scales, \( \varepsilon(h) \) and \( \varepsilon^2(h) \). More precisely, they assume the existence of the limits
\begin{equation} \label{eq:098777}
\gamma_1 := \lim_{h \to 0} \frac{h}{\varepsilon(h)} \quad \text{and} \quad \gamma_2 := \lim_{h \to 0} \frac{h}{\varepsilon^2(h)},
\end{equation}
and analyze three cases:
\begin{itemize}
    \item \( \gamma_1 = 0 \), \( \gamma_2 = +\infty \),
    \item \( \gamma_1 \in (0, +\infty) \), \( \gamma_2 = +\infty \),
    \item \( \gamma_1 = \gamma_2 = +\infty \).
\end{itemize}
The identification of the reduced models for $\gamma_1 = 0$ and $\gamma_2 \in [0,+\infty)$ remains an open problem in \cite{BUFDAVFON15}. In \cite{VEL15}, the critical case $h\sim \varepsilon^2(h)$, i.e. $\gamma_1=0$, $\gamma_2\in(0,+\infty)$ is studied. The more delicate case $\gamma_1=\gamma_2=0 $ has only been partially treated in \cite{HORVEL15} for shells, but only in the one-scale von Kármán regime.\\
In this paper, we extend this framework to the bending shell theory by including multiscale effects.
We derive different asymptotic theories depending on the parameters \( \gamma_1 \) and \( \gamma_2 \), as defined in \eqref{eq:098777}, and consider the following cases:
\begin{itemize}
    \item[(a)] \( \gamma_1 \in [0, +\infty] \), \( \gamma_2 = +\infty \): includes the subcases \( \gamma_1 = 0 \), \( \gamma_1 \in (0, +\infty) \), and \( \gamma_1 = \gamma_2 = +\infty \).
    \item[(b)] \( \gamma_1 = 0 \), \( \gamma_2 \in (0, +\infty) \).
    \item[(c)] \( \gamma_1 = \gamma_2 = 0 \).
\end{itemize}
We prove case $(a)$ using standard arguments (see \cite{HORVEL15,HORVEL18}), while cases $(b)$ and $(c)$ are addressed through a novel approach.
More precisely, we consider a shell \( S^h \subset \mathbb{R}^3 \) of thickness \( h \in (0,1 ]\), subjected to two material heterogeneities scales \( \varepsilon(h) \) and \( \varepsilon^2(h) \), according to the assumption \eqref{eq:098777}. Using a rescaling \( \Theta^h : S^h \rightarrow S^1 \) and a chart inverse \( r \) defining the shell's mid-surface \( S \), we define the rescaled nonlinear elastic energy as
\begin{equation}\label{eq1:juni03}
J^h(v) = \frac{1}{h} \int_{S^h} W\left(\Theta^h(x), \frac{r(x)}{\varepsilon}, \frac{r(x)}{\varepsilon^2}, \nabla v(x)\right) \, dx,
\end{equation}
for every deformation \( v \in W^{1,2}(S^h; \mathbb{R}^3) \), where the stored energy density \( W \) is periodic in its second and third arguments and satisfies the standard assumptions of nonlinear elasticity, including a non-degeneracy condition (see Section~\ref{sec:4}), as in \cite{HORVEL18, HORVEL15}. Let \( S^1 \) be the shell of thickness 1, and set \( Y = [-\tfrac{1}{2}, \tfrac{1}{2})^2 \). Note that here, in the definition of energy $J^h$ (\ref{eq1:juni03}), we write as $r(x)$ by the simplicity of notation, but it needs to be interpreted as $r(\pi(x))$ to be able to map the shell of thickness one $S^1$ to the mid surface of the plate $\omega$ (see Remark \ref{obs:may12}), where $\pi$ is the nearest of retraction point function which maps $S^h$ (with $h$ small) to its midsurface $S$ (see Definition \ref{nearestpointofretrac}). We denote by \( \mathcal{Y} = \mathcal{Z} = \mathbb{R}^2 / \mathbb{Z}^2 \) the torus associated with \( Y \). For each \( (x, y, z) \in S^1 \times \mathcal{Y} \times \mathcal{Y} \), we assume the existence of a quadratic form 
\[
\mathscr{Q}(x, y, z, \cdot): \mathbb{R}^{3 \times 3} \to \mathbb{R}
\]
representing the linearization of \( W \), such that
\begin{equation} \label{eq:linW}
\lim_{\|\mathbf{G}\| \to 0} \, \mathrm{ess\,sup}_{(x,y,z) \in S^1 \times \mathcal{Y} \times \mathcal{Y}} \frac{|W(x,y,z, \mathrm{Id} + \mathbf{G}) - \mathscr{Q}(x,y,z, \mathbf{G})|}{\|\mathbf{G}\|^2} = 0.
\end{equation} For \( \gamma_1 \in [0, +\infty] \), \( \gamma_2 \in [0, +\infty] \), and \( x \in S \), we define the functional
\[
\mathscr{Q}_{\gamma_1, \gamma_2}(x, q) = \inf \int_I \int_{\mathcal{Y}} \int_{\mathcal{Y}} \mathscr{Q}\left(x + t n(x), y, z, p + t q + U(t, y, z)\right) \, dz \, dy \, dt,
\]
where the infimum is taken over all 
\[
U \in L^{(x)}_{\gamma_1, \gamma_2}(I \times \mathcal{Y} \times \mathcal{Y})
\]
(see Section~\ref{subsectiongammas} for the precise definitions) and all \( p \in T^*_x S \otimes T^*_x S \), with \( I = \left(-\tfrac{1}{2}, \tfrac{1}{2}\right) \).
We then define the limit energy functional
\[
I_{\gamma_1, +\infty}(u) = 
\begin{cases}
\int_S \mathscr{Q}_{\gamma_1, +\infty}(x, \mathbf{S}_u^r(x)) \, d\mathrm{vol}_S(x) & \text{if } u \in W^{2,2}_{\mathrm{iso}}(S), \\
+\infty & \text{otherwise}.
\end{cases} \]

To improve readability, we present our main result here (for details, see Section \ref{sec:4}, Theorem \ref{teo:187}, and Section \ref{sec:6}, Theorem \ref{theo61}):\\
\begin{teo}
Let  $(u^h)\subset W^{1,2}(S^h;\RR^3)$ be a sequence of deformations satisfying
\begin{equation}\label{eq:33}
\displaystyle\lim_{h\to 0}\sup h^{-2}J^h(u^h)<+\infty.
\end{equation}
Then 
\begin{equation}\label{mainresultmay10}
h^{-2} J^h(u^h)\quad\Gamma-\text{converges to}\quad I_{\gamma_1, \gamma_2}(u)\,\text{as}\, h\to 0,
\end{equation}
in the sense of \cite{REF6},
i.e. 
\begin{itemize}
\item[(i)](\textbf{Liminf inequality}) Let $(u^h)\subset W^{1,2}(S^h,\RR^3)$ be such that (\ref{eq:33}) holds and such that $y^h-\frac{1}{|S^1|}\int_{S^1}y^h\to u$ strongly in $L^2(S^1)$ for some $u\in L^2(S^1, \RR^3)$. Then
\begin{equation*}
\lim_{h\to 0}\inf h^{-2} J^h(u^h)\geq I_{\gamma_1, \gamma_2}(u).
\end{equation*}
\item[(ii)](\textbf{Recovery sequence}) If, in addition, $S$ is simply connected, then for every $u\in \widetilde{W}^{2, 2}_{\text{iso}}(S)$ there exists $(u^h)\subset W^{1,2}(S^h;\RR^3)$ satisfying (\ref{eq:33}), and such that $y^h\to u$ strongly in $W^{1,2}(S^1)$. Moreover,
\begin{equation*}
\lim_{h\to 0}h^{-2}J^h(u^h)=I_{\gamma_1,\gamma_2}(u).
\end{equation*}
\end{itemize}
\end{teo}
The key to obtaining this $\Gamma-$convergence result remains in identifying the effective energy associated with the rescaled elastic energies
\begin{equation}\label{energyaim:may12}
\displaystyle\frac{J^h(u^h)}{h^2}=\displaystyle\frac{1}{h^2|S^h|}\int_{S^h}W(x,r(x)/\varepsilon, r(x)/\varepsilon^2, \nabla u^h(x))dx,
\end{equation}
for different values of $\gamma_1$ and $\gamma_2$. 
Our analysis relies on the following approaches:
\begin{itemize}
\item[(1)] Dimension reduction techniques, in particular the quantitative rigidity estimate and approximation schemes developed by Friesecke, James, and Müller in their work on the derivation of nonlinear plate theories \cite{FJM02}.
\item[(2)] Homogenization methods, especially three-scale convergence (see \cite{Allamults, BUFDAVFON15, ferreiramsc}, and \cite{fonsecamsc}).
\end{itemize}

Starting from a limiting deformation in the space $W^{2,2}(S^h; \mathbb{R}^3)$, we construct the required recovery sequence using the density result from \cite{HORVEL18}. One of the main challenges lies in addressing the technical complexities that arise across all scaling regimes, due to the interaction of three distinct scales. Each case, characterized by specific values of $\gamma_1$ and $\gamma_2$, is treated using tailored techniques. Some cases are handled with standard tools; in others, the shell's curvature introduces new terms absent in the plate setting. In yet other situations, we restrict our attention to particular classes of shells to manage the complexity. Below, we outline the methods used in each case to overcome these difficulties.\\

The most important task is to identify the three-scale limit of the sequence of linearized stresses $\widetilde{\textbf{G}^h}$, given by
\begin{equation}\label{eq:34}
\widetilde{\textbf{G}^h} = \displaystyle\frac{(R^h)^T\nabla_h y^h - I}{h},
\end{equation}
where $R^h$ is a sequence of rotations, $I$ is the $3 \times 3$ identity matrix, and $y^h$ is the rescaled version of the deformation $u^h$ (see Section \ref{subsect51} for more details).
The three-scale analysis improves the technical understanding of the problem across all scaling regimes.
The novelty of this work, compared to the plate case (see \cite{BUFDAVFON15} and \cite{VEL15}), is that furthermore of be able to deal the standard cases $\gamma_1\in[0,+\infty]$ y $\gamma_2=+\infty$, we can address the cases $\gamma_1 = 0$, $\gamma_2 \in [0, +\infty)$,
i.e., the cases $\gamma_1 = 0$, $\gamma_2 \in (0, +\infty)$, and $\gamma_1 = \gamma_2 = 0$. More precisely,
\begin{itemize}
    \item  For $\gamma_1 \in (0,+\infty]$ and $\gamma_2 = +\infty$, we use the Friesecke–James–Müller rigidity estimate (see \cite{FJM02}, Theorem 4.1), which provides sequences of rotations $\{R^h\}_h$ that are piecewise constant on cubes of size $\varepsilon(h)$ centered at points in $\varepsilon(h)\mathbb{Z}^2$.
To identify the three-scale limit of the linearized stresses, it is necessary to consider oscillating sequences on the finer scale $\varepsilon^2(h)$.
We address this by subdividing the cubes of size $\varepsilon^2(h)$, with centers in $\varepsilon^2(h)\mathbb{Z}^2$, into so-called “good cubes” and “bad cubes” (see \cite{BUFDAVFON15}, Theorem 4.2.1). In fact, good cubes are those that remain completely contained within a larger cube of size $\varepsilon(h)$ centered at a point in $\varepsilon(h)\mathbb{Z}^2$.
We show that the measure of the intersection between the domain $\Omega$ and the set of bad cubes tends to zero faster than $\varepsilon(h)$ as $h \to 0$. This result allows us to establish a lower bound for the rescaled nonlinear elastic energy $J^h$.
\item For  $\gamma_1=0$ and $\gamma_2=+\infty$, again, the Friesecke–James–Müller rigidity rigidity estimate (see  \cite{FJM02}, Theorem 4.1) provides sequences of rotations  $\{R^h\}_h$ which are piecewise constants on cubes of sizes $\delta = \displaystyle\left(2\lceil\frac{h}{\varepsilon^2(h)}\rceil+1\right)\varepsilon^2(h)$, where $\lceil\cdot \rceil$ is the integer part. 
To define new rotations that are piecewise constant on cubes of size $\varepsilon^2(h)$, we combine these with the results from \cite{BUFDAVFON15} and \cite{HORVEL18} (included in the Appendix for the reader's convenience).
These new rotations allow us to obtain the multiscale limit of the stress, which then yields a lower bound for the sequence of stored energy functionals.
Interestingly, in this case—unlike the previous one—the identification of good and bad cubes of size $\varepsilon^2(h)$ is not useful, as the oscillations of the test functions over the cubes are not negligible.
\item For the open case $\gamma_1 = \gamma_2 = 0$, i.e., $h \ll \varepsilon^2$ (see Section \ref{sec:6}), which includes multiscale effects, we restrict our analysis to convex shells, thus highlighting the stronger role of geometry in this case.
We use the Friesecke–James–Müller rigidity estimate to construct rotations that are piecewise constant on cubes of size $\delta(h) = h$. However, here we also need to subdivide those cubes into “good cubes” (which are fully contained in a cube of size $\varepsilon^2$ belonging to $Q_{\varepsilon^2}$), where
\begin{equation*}
Q_{\varepsilon^2}:=\displaystyle\bigcup_{z\in \ZZ^{\varepsilon^2}}Q(\varepsilon^2(h) z, \varepsilon^2(h)).
\end{equation*}
and “bad cubes” (those intersecting more than one element of $Q_{\varepsilon^2}$).
Finally, we complete the characterization of the three-scale limit of the linearized stresses by arguing as in the case $\gamma_1 \in (0, +\infty)$ and $\gamma_2 = +\infty$.
\item For $\gamma_1=0$ and $\gamma_2\in (0,+\infty)$ since the Weingarten map $\textbf{S}$ is non-zero, due to the shell structure, this activates the term $\gamma_2^{-1}\varphi\,\textbf{S}(x) $
which vanishes in the plate case (where $\textbf{S} = 0$), see Remark \ref{remarknewmay12}.
We use the Friesecke–James–Müller rigidity estimate to find rotations that are piecewise constant on cubes of size $\delta(h) = \varepsilon^2(h)$, and we proceed similarly to the case $\gamma_1 = 0$ and $\gamma_2 = +\infty$.
\end{itemize}
In comparison to Refs. \cite{Neukammvelcic} and \cite{VEL15}, we deal with the structure of the homogenized energy which is obtained by means of a double pointwise minimization, first with respect to the faster periodicity scale, and then with respect to the slower one and the $x_3$ variable density. Mathematically speaking, we use three-scale techniques in place of two-scale. In homogenization theory, two-scale and three-scale convergence rigorously analyze PDEs with multi-scale oscillations, common in composites and porous media (see for instance \cite{BUZANCICHERNANDEZ2025}). Two-scale convergence \cite{N89,REF0} links macroscopic variable $x $ with a single microscopic scale $y = x/\varepsilon $, ideal for periodic structures. However, many systems (e.g., hierarchical porous media) require additional scales due to nested microstructures. Three-scale convergence \cite{Allamults} introduces a second fast variable $z = x/\varepsilon^2 $, capturing effects like heat conduction in materials with pores at different scales. While two-scale uses one microscopic variable, three-scale extends this to multiple nested scales, enabling precise modelling of complex heterogeneous media.

\vspace{0.25 cm}

This work is relevant because it plays a key role in the study of complex nonlinear elastic materials that are microscopically heterogeneous, such as composites, foams, shape-memory alloys, or materials with periodic microstructures. This type of model enables a rigorous and efficient connection between the material's microscopic behavior (microscale) and its macroscopic response (macroscale). Moreover, it is applicable to a wide range of scenarios, including:

\begin{itemize}
\item Biomechanics: Tissues such as bone or skin exhibit hierarchical microstructures across multiple scales (see, for instance, \cite{Millerpenta2025}).
\item Civil engineering: Materials such as concrete or geotechnical composites contain heterogeneities spanning several orders of magnitude.
\end{itemize}

\vspace{0.5 cm}

The paper is organized as follows: Section \ref{sec:2} introduces the definition and some properties of three-scale convergence for shells. Section \ref{sec:3} describes the setting of our problem and presents the assumptions on the energy density. In Section \ref{sec:4}, we state our main result, which is proved in Section \ref{sec:5}, where we derive bounds for the rescaled elastic energy depending on different values of $\gamma_1$ and $\gamma_2$. Finally, Section \ref{sec:6} analyzes the homogenized model in the case of convex shells with $h \ll \varepsilon^2$.
\section{Preliminaries}\label{sec:2}
\subsection{Homogenization for shells}
We begin by introducing some further notation. We use $\mathcal{Y}$ and $\mathcal{Z}$ to refer us  to the scales $y=x/\varepsilon$ and $z=x/\varepsilon^2$ respectively. For all $k\in\NN\cup \{0\}$ the set of all $f\in C^{k}(\RR^2)$ with $D^{\alpha}f(\cdot+z)=D^{\alpha}f$ for all $z\in\ZZ^2$ and all multiindices $\alpha$ of order up to $k$ is denoted by $C^{k}(\mathcal{Y})$. For $h\in(0,1]$, let $m(h)$ and $n(h)$ two positive functions from $(0,1]$ to $(0,\infty)$. We denote by $Q(m(h)z,n(h))$, the square parallelepiped grids with a square base of each side $m(h)$ and height $n(h)$ where $z\in\ZZ^2$.

We denote by $C^k_0$ functions with compact support. For any open set $A$, we denote by $L^2(\mathcal{Y})$, $W^{1,2}(\mathcal{Y})$ and $W^{1,2}(A\times\mathcal{Y})$ the Banach spaces obtained as closures of $C^{\infty}(\mathcal{Y})$ and $C^{\infty}(\overline{A}, C^{\infty}(\mathcal{Y}))$ with respect to the norm in $L^2(Y)$, $W^{1,2}(Y)$ and $W^{1,2}(A\times Y)$, respectively. An additional dot (e.g. in $\dot{L}^2(\mathcal{Y})$) denotes functions with average zero over $Y$.
\subsubsection{Surfaces and shells in $\RR^3$}
In this section, we introduce the notations and terminologies related to domain and tools necessary for the analysis of this current article. In fact, we have taken them from \cite{HORVEL15,HORVEL18} but for the convenience of readers, we reproduce it here.

Let $h\in (0,1]$, $\kappa\in(0,1)$ and let $\omega\subset\RR^2$ be a bounded domain with $C^{3, \kappa}$ boundary. Set $I=(-\frac{1}{2}, \frac{1}{2})$, $\Omega^h=\omega\times (hI)$, and $\Omega=\omega\times I$. From now on, $S\subset\RR^3$ denotes (the relative interior of) an embedded compact connected oriented surface with boundary. For convenience, we assume that a single chart parametrizes $S$. More precisely, we assume that there exists an open set $V\subset \RR^3$ containing the closure of $S$ and an open set $U\subset\RR^3$ containing $\overline{\omega}\times\{0\}$ and $C^{3,\kappa}$ diffeomorphism $\Phi: V\to U$ such that.
\begin{equation*}
\Phi(S)=\omega\times\{0\}.
\end{equation*}

Then $\xi:\omega\to S$, defined by $\xi(z)=\Phi^{-1}(z,0)$, is a global $C^{3,\kappa}$ chart for $S$.\\
By $W^{2,2}_{\text{iso}}(S)$ we denote the $W^{2,2}(S)$ isometries of the surface $S$ into $\RR^3$. The space $W^{2,\infty}_{\text{iso}}(S)$ is defined similarly. Clearly $u\in W^{2,2}_{\text{iso}}(S)$ is equivalent to $u\circ \xi\in W^{2,2}_{g}(\omega)$, for $g=(\nabla \xi)^{T}(\nabla\xi)$ the Riemannian metric on $\omega$ induced by $\xi$.

\medskip As usual, $T\,S$ denotes the tangent bundle over $S$ and $N\, S$ the normal bundle. A basis of the tangent space $T_{x}S$ is given by
\begin{equation*}
\tau_{\beta}(x)=(\partial_{\beta}\xi)(\Phi(x))\quad\text{for all}\,x\in S,
\end{equation*}
where $\beta=1, 2$. We view $T_{x}S$ as a subspace of $\RR^3$ and write $\sigma\cdot\tau$, the scalar product of both spaces.

\medskip The dual basis of the tangent space $T_{x}S$ is denoted by $(\tau^1(x), \tau^2(x))$. So by definition
\begin{equation*}
\tau^{\alpha}\cdot\tau_{\beta}=\delta_{\alpha\beta}\,\text{on}\, S,
\end{equation*}
where $\delta_{\alpha\beta}$ is the Kronecker delta symbol. We frenquently identify $T^{*}_{x}(S)$ with $T_{x}(S)$ via the scalar product.
\medskip Define the normal $n: S\to \mathbb{S}^2$ by
\begin{equation*}
n=\displaystyle\frac{\tau_1\times\tau_2}{|\tau_1\times\tau_2|}.
\end{equation*}
The orthogonal projection onto $T_{x}S$ is
\begin{equation*}
T_{S}(x)=I-n(x)\otimes n(x).
\end{equation*}
The tensor products $T\,S\otimes T\, S$ etc. are defined fiberwise. $T^{*}_{x}S\otimes T^{*}_{x}S$ will be regarded as a subspace of $\RR^{3\times 3}$.\\
If $E$ and $F$ are vector spaces (or bundles) then the space of all symmetric products
\begin{equation*}
a\odot b:=\displaystyle\frac{1}{2}(a\otimes b+b\otimes a),
\end{equation*}
with $a\in E$ and $b\in F$ is denoted by $E\odot F$.

Sections $B$ of $T^{*}S\otimes T^{*} S$ will frequently regarded as maps from $S$ into $\RR^{3\times 3}$ via the embedding $\iota$ defined by $\iota(B)=B(T_{S}, T_{S})$. By definition, $B(T_{S}, T_{S}) : S\to \RR^3$ takes the vector fields $v, w : S\to \RR^3$ into the function $x\mapsto B(x)(T_{S}(x) v(x), T_{S}(x) w(x))$.

For any vector bundle $E$ over $S$ we denote by $L^2(S, E)$ the space of all $L^2-$sections of $E$. The spaces $W^{1,2}(S, E)$ etc. are defined similarly. For any vector bundle $E$ over $S$ with fibers $E_{x}$, we denote by $L^{2}(\mathcal{Y}, E)$ the vector bundle over $S$ with fibers $L^2(\mathcal{Y}, E_{x})$. The bundles $W^{1,2}(\mathcal{Y}, E)$ etc. are defined similarly. For example, $L^2-$sections of the bundle $W^{1,2}(\mathcal{Y}, T\,S)$ are given by
\begin{align*}
L^2(S,& W^{1,2}(\mathcal{Y}, T\,S))=\\
&\{Z\in L^2(S, W^{1,2}(\mathcal{Y}, \RR^3)): Z(x)\in W^{1,2}(\mathcal{Y}, T_{x}S)\,\,\text{for a.e.}\,\, x\in S \}.
\end{align*}

For a function $f : S\to \RR$ its differential $df$ is given by $df(x)\tau=\nabla_{\tau}f(x)$ for all $\tau\in T_{x}S$. Here $\nabla_{\tau}f$ denotes the directional derivative of $f$ in the direction of the tangent vector $\tau$. We extend these definitions componentwise to maps into $\RR^3$. By $\nabla$ we denote the usual gradient on $\RR^3$ (or on $\RR^2$).

As usual, the Weingarten map $\textbf{S}$ of $S$ is the differential of the normal, \textit{i.e.},
\begin{equation*}
\textbf{S}(x)\tau=(\nabla_{\tau}n)(x)\,\,\text{for all}\,\, x\in S,\,\tau\in T_{x}S.
\end{equation*}
We extend $\textbf{S}(x)$ trivially to $\RR^3$ by setting $\textbf{S}(x)=\textbf{S}(x)T_{S}(x)$.

For an immersion $u : S\to \RR^3$ denote by $\textbf{S}_{u}$ the Weingarten map for the surface $u(S)$. We define its pullback to $S$ by setting
\begin{equation*}
(u^{*}\textbf{S}_{u})\tau=u^{*}(\textbf{S}_{u}D_{\tau}u)
\end{equation*}
for all smooth tangent vector fields $\tau$ to $S$. Here by definition, $u^{*}(D_{\sigma}u)=\sigma$ for all smooth tangent vector fields $\sigma$ to $S$. Following the ideas given in \cite{FJM03}, we define the relative Weingarten map
\begin{equation*}
\textbf{S}^{r}_{u}=u^{*}\textbf{S}_{u}-\textbf{S}.
\end{equation*}

Using the Weingarten map, the covariant derivative of a tangent vector field $\tau$  along another tangent vector field
$\sigma$ is given by
\begin{equation*}
D_{\sigma}\tau=\nabla_{\sigma}\tau+\langle \textbf{S}\tau, \sigma\rangle n,
\end{equation*}
or simply: $D\tau=\nabla \tau+n\otimes \textbf{S}\tau$.

With a given displacement $V : S\to \RR^3$ one associate the following quantity:
\begin{itemize}
\item The quadratic form $q_V$ given by
\begin{equation*}
q_V(\tau,\eta)=\displaystyle\frac{1}{2}(\eta\cdot \nabla_{\tau}V+\tau\cdot \nabla_{\eta}V),
\end{equation*}
for all tangent vector fields $\tau,\eta$ along $S$.
\item For tangent vector fields $v$ along $S$ the quadratic form corresponding to $q_v$ is denoted
by $\text{Def}_S v$ and called deformation tensor of $v$. It is given by the Lie-derivative of the metric in direction $v$, i.e.,
\begin{equation*}
(\text{Def}_S v)=\displaystyle\frac{1}{2}\left(\eta\cdot D_{\tau}v+\tau\cdot D_{\eta}v\right)=\displaystyle\frac{1}{2}\left(\eta\cdot \nabla_{\tau}v+\tau\cdot \nabla_{\eta}v\right),
\end{equation*}
for all tangent vector fields $\tau$ and $\sigma$.
\end{itemize}

The following result is essentially contained in Lemma 2.1 of \cite{HORVEL15}.
\begin{lemma}
Let $V\in W^{1,2}(S;\RR^3)$. The we have almost everywhere on $S$,
\begin{equation*}
q_V=\text{Def}_S V_{tan}+(V\cdot n)\textbf{S}.
\end{equation*}
\end{lemma}
 
 \begin{definition}\label{nearestpointofretrac}
We denote by $\pi$ the nearest point of retraction $\pi$ of a tubular neighborhood of $S$ onto $S$ which satisfies $\pi(x+tn(x))=x$ for small $|t|$ and all $x\in S$.     
 \end{definition}
 
 After rescaling the ambient space, we may assume that the curvature of $S$ is as small as we please. Therefore, we may assume without loss of generality that $\pi$ is well-defined on a domain containing the closure of the set $\{x+tn(x): x\in S, -1/2<t<1/2\}$, and that $|Id+t\textbf{S}(x)|\in (1/2, 3/2)$ for all $t\in(-\frac{1}{2}, \frac{1}{2})$ and all $x\in S$.

 For a subset $S_0\subset S$ and $h\in (0,1]$ we define $S^h_0=\{x+tn(x): x\in S_0,\, -h/2<t<h/2\}$. In particular, the whole shell is, by definition,
 \begin{equation*}
 S^h=\left\{x+tn(x): x\in S\,\text{and}\, t\in\displaystyle(-\tfrac{h}{2}, \tfrac{h}{2})\right\}.
 \end{equation*}
We introduce the map $r=\Phi\circ \pi$. Moreover, we introduce the function $t: S^1\to\RR$ by setting $t(x)=(x-\pi(x))\cdot n(x)$ for all $x\in S^1$. We have the following identity on $S^1$, cf.  \cite{HORVEL15}
\begin{equation}\label{eq:187}
d\pi=T_{S}(\pi)(I+t\textbf{S}(\pi))(I+t\textbf{S}(\pi)T_{S}(\pi))^{-1}.
\end{equation}
(Here and elsewhere we write $T_{S}(\pi)$ instead of $T_{S}\circ \pi$ etc.). Hence there exists a constant $C$ depending only on $S$ such that
\begin{equation}\label{eq:281}
|d\pi-(I-t\textbf{S}(\pi)T_{S}(\pi))|\leq C t^2\,\,\text{on}\,\, S^1.
\end{equation}

\begin{obs}\label{obs:may12}
Abusing notations, maps $f: S\to \RR^{k}$ will often be extended to $S^1$ by setting $f=f\circ \pi$. We extend $r, T_{S}$ and $\textbf{S}$ in this way, too.  
\end{obs}

 For functions $f\in L^2(S, W^{2,2}(\mathcal{Y}))$ the expression $\text{Hess}_{\mathcal{Y}}$ is the section of the bundle $L^2(\mathcal{Y}, T\, S\odot T\,S)$ over $S$ given by
 \begin{equation*}
(\text{Hess}_{\mathcal{Y}}f)(x,y)=(\nabla^2_yf)_{\alpha\beta}(x,y)\tau^{\alpha}(x)\odot\tau^{\beta}(x),
\end{equation*}
where $(\nabla^2_{y}f)_{\alpha\beta}=\partial_{y_{\alpha}}\partial_{y_{\beta}}f$. Analogously,  for functions $f\in L^2(S\times\mathcal{Y}, W^{2,2}(\mathcal{Y}))$ the expression $\text{Hess}_{\mathcal{Z}}$ is the section of the bundle $L^2(\mathcal{Y}\times \mathcal{Y}, T\, S\odot T\,S)$ over $S$ given by
 \begin{equation*}
(\text{Hess}_{\mathcal{Z}}f)(x,y, z)=(\nabla^2_zf)_{\alpha\beta}(x,y,z)\tau^{\alpha}(x)\odot\tau^{\beta}(x),
\end{equation*}
where $(\nabla^2_{z}f)_{\alpha\beta}=\partial_{z_{\alpha}}\partial_{z_{\beta}}f$.

 \medskip For $v\in L^2(S, W^{1,2}(\mathcal{Y};\RR^2))$ we define the section $\text{Def}_{\mathcal{Y}} v$ of the bundle $L^2(\mathcal{Y}, T^{*}S\odot T^{*}S)$ by
\begin{equation*}
(\text{Def}_{\mathcal{Y}}v)(x,y)=\left(\text{sym}\,\nabla_{y}v(x,y)\right)_{\alpha\beta}\tau^{\alpha}(x)\odot\tau^{\beta}(x).
\end{equation*}
For $v\in L^2(S\times\mathcal{Y}, W^{1,2}(\mathcal{Y};\RR^2))$ we define the section $\text{Def}_{\mathcal{Z}} v$ of the bundle $L^2(\mathcal{Y}\times\mathcal{Y}, T^{*}S\odot T^{*}S)$ by
\begin{equation*}
(\text{Def}_{\mathcal{Z}}v)(x,y, z)=\left(\text{sym}\,\nabla_{z}v(x,y, z)\right)_{\alpha\beta}\tau^{\alpha}(x)\odot\tau^{\beta}(x).
\end{equation*}
Here and elsewhere $\nabla_{y}$ is the gradient in $\mathcal{Y}$ with respect to the variable $y$ and $\nabla_{z}$ is the gradient in $\mathcal{Z}$ with respect to the variable $z$ (and not some directional derivative).

 \medskip We define the map $\Xi :\omega\times \RR\to\RR^3$ by
 \begin{equation*}
 \Xi(z', z_3)=\xi(z')+z_3 n\left(\xi(z')\right)\,\,\text{for all}\,\,z'\in\omega\,\,\text{and}\,\,z_3\in\RR.
 \end{equation*}
 
 As in \cite{HORVEL14}, we will use the diffeomorphism $\widetilde{\Theta^h}:\Omega^h\to \Omega$ given by $\widetilde{\Theta^h}(z_1, z_2, z_3)=(z_1, z_2, z_3/h)$, and for a map $\widetilde{y}:\Omega\to\RR^3$ we introduce the scaled gradient $\widetilde{\nabla}_hy=(\partial_1 y, \partial_2 y, \frac{1}{h}\partial_3 y)$. The counterpart of $\widetilde{\Theta^h}$ on the shell is the diffeomorphism $\Theta^h:S^h\to S^1$ given by
 \begin{equation*}
 \Theta^h(x)=\pi(x)+\displaystyle\frac{t(x)}{h}n(x).
 \end{equation*}
 It is easy to see that
 \begin{equation*}
 \Theta^h\circ\Xi=\Xi\circ\widetilde{\Theta}^h\quad\text{on}\,\Omega^h.
 \end{equation*}
 For given $u: S^h\to\RR^3$ we define its pulled back version $\widetilde{u}:\Omega^h\to \RR^3$ by $\widetilde{u}=u\circ\Xi$. We also define its rescaled version $y: S^1\to\RR^3$ by $y(\Theta^h)=u$ on $S^h$ and we define the pulled back version $\widetilde{y}$ of this map by $\widetilde{y}=y\circ\Xi$. Then it is easy to see that
 \begin{equation}\label{eq:2100087}
(\widetilde{\nabla}_h\widetilde{y})\circ\widetilde{\Theta}^h=\nabla\widetilde{u}\quad\text{on}\,\Omega.
 \end{equation}
 We define the rescaled gradient $\nabla_{h}y$ of $y$ by the condition
 \begin{equation}\label{eq:2200087}
 (\nabla_h y)\circ \Theta^h=\nabla u\quad\text{on}\,S^h.
 \end{equation}
 Using (\ref{eq:2100087}) and (\ref{eq:2200087}) it is clear to see that
 \begin{equation}\label{eq:2300087}
 \widetilde{\nabla}_h\widetilde{y}=\nabla u(\Xi)\left((\nabla\Xi)\circ(\widetilde{\Theta}^h)^{-1}\right)
 \end{equation}
 and using (\ref{eq:187}) we can see that
\begin{equation}\label{eq:387}
\nabla \Theta^h=\left(T_{S}+\displaystyle\frac{1}{h}(n\otimes n+t\textbf{S})\right)(I+t\textbf{S})^{-1}\,\,\text{on}\,\, S^h.
\end{equation}
Finally, to express $\nabla_h y$ in terms of $\nabla y$, we insert the definition of $y$ into (\ref{eq:2200087}) and use (\ref{eq:187}) to find
\begin{equation}\label{eq:2500087}
\nabla_hy=\nabla y\left(T_{S}+\displaystyle\frac{1}{h}n\otimes n+t\textbf{S}\right)(I+ht\textbf{S})^{-1}\quad\text{on}\, S^1.
\end{equation}

\subsubsection{Three-scale convergence on shells}

Taking Velčić's two-scale convergence given in \cite{HORVEL15} as motivation, we can define a version for three-scale convergence in curved domains. Recall that we extend the chart $r$ trivially from $S$ to $S^1$. We make the following definitions:

\begin{itemize}

\item[(i)] A sequence $(f^h)\subset L^2(S^1)$ is said to converge weakly three-scale on $S^1$ to the function $f\in L^2(S^1, L^2(\mathcal{Y}\times \mathcal{Y}))$ as $h\to 0$, provided that the sequence $(f^h)$ is bounded in $L^2(S^1)$ and
\begin{equation}\label{eq:1****}
\displaystyle\lim_{h\to 0}\int_{S^1}f^h(x)\varphi\left(x,\frac{r(x)}{\varepsilon},\frac{r(x)}{\varepsilon^2}\right)dx=\int_{S^1}\int_{\mathcal{Y}}\int_{\mathcal{Y}}f(x, y, z)\varphi\left(x, y, z\right)\,dz\,dy\,dx,
\end{equation}

for all $\varphi\in C^{0}_{c}\left(S^1,C^0(\mathcal{Y}\times \mathcal{Y}) \right)$.

\item[(ii)] We say that $f^h$ strongly three-scale converges to $f$ if, in addition,
\begin{equation*}
\displaystyle\lim_{h\to 0}||f^h||_{L^2(S^1)}=||f||_{L^2(S^1\times \mathcal{Y}\times \mathcal{Y})}.
\end{equation*}

\item[(iii)] For a sequence $(f^h)\subset L^2(S^1)$ and for $f_1\in L^2(S^1\times \mathcal{Y}\times \mathcal{Y})$ with $\int_{\mathcal{Y}}f_1(\cdot, \cdot,z)dz=0$ for almost every $x\in S^1\times \mathcal{Y}$, we write $f^h\overset{osc, Z}{\rightharpoonup} f_1$ provided that
\begin{equation}\label{eq:2****}
\displaystyle\int_{S^1}f^h(x)\varphi\left(x, r(x)/\varepsilon\right)\rho\left(r(x)/\varepsilon^2\right)dx\to\int_{S_1}\int_{\mathcal{Y}}\int_{\mathcal{Y}}f_1(x,y,z)\varphi(x,y)\rho(z)\,dz\,dy\,dx
\end{equation}
\end{itemize}
for all $\varphi\in C^{\infty}_0(S^1; C^{\infty}(\mathcal{Y}))$, all $\rho\in C^{\infty}(\mathcal{Y})$  with $\int_{\mathcal{Y}}\rho(z) dz=0$.

\medskip We write $f^h\overset{3}{\rightharpoonup} f$ to denote weak three-scale convergence and $f^h\overset{3}{\rightarrow} f$ to denote strong three-scale convergence. If $f^h\overset{3}{\rightharpoonup} f$ then $f^h\rightharpoonup \int_{\mathcal{Y}}\int_{\mathcal{Y}}f(\cdot, y, z)dzdy$ weakly in $L^2$. If $f^h$ is bounded in $L^2(S^1)$, then it has a subsequence which converges weakly three-scale to some $f\in L^2(S^1; L^2(\mathcal{Y}\times \mathcal{Y}))$. Other results can be inferred from their counterparts in planar domains (cf. \cite{REF0} and \cite{REF33}) using the following basic observations.

\vspace{0.5 cm}
\begin{obs}\label{remarkthreescaleconvmarch21}
\medskip Defining $\tilde{f}=f^h\circ \Xi$ and $\tilde{f}(z,y,w)=f(\Xi(z), y, w)$, and taking
\begin{equation*}
\tilde{\varphi}(z, y, w)=\varphi\left(\Xi(z), y, w\right)\left(\det \nabla\Xi^T(z)\nabla\Xi(z)\right)^{1/2},
\end{equation*}
a change of variable shows that (\ref{eq:1****}) is equivalent to
\begin{equation}\label{eq:3****}
\displaystyle\int_{\Omega}\tilde{f}^h(z)\tilde{\varphi}\left(z, z'/\varepsilon, z'/\varepsilon^2\right)dz\to\int_{\Omega}\int_{\mathcal{Y}}\int_{\mathcal{Y}}\tilde{f}(z,y,w)\tilde{\varphi}(z,y,w)\,dw\,dy\,dz,
\end{equation}
where $z'$ is the projection of $z$ onto $\RR^2$. Therefore $f^h\overset{3}{\rightharpoonup} f$ on $S^1$ if and only if $\tilde{f}^h\overset{3}{\rightharpoonup} \tilde{f}$ on $\Omega$ in the usual sense.
\end{obs}
\vspace{0.25 cm}

\begin{obs}\label{oscilationandthrescale}
As a direct consequence of the definition of multiscale convergence and density arguments, if $\{\tilde{f}^h\}\subset L^2(\omega)$, then
\begin{equation*}
\tilde{f}^h\overset{3}{\rightharpoonup}\tilde{f}\,\,\text{weakly $3-$scale in}\,\, L^2(\Omega)
\end{equation*}
if and only if
\begin{equation*}
f^h(x)\overset{osc, Z}{\rightharpoonup}\tilde{f}-\displaystyle\int_{Y}\tilde{f}(x,y,z)\,dz.
\end{equation*}
\end{obs}

\vspace{0.25 cm}

The following lemma can be also adapted to the curved setting following the same pattern as above (writing just $S^1$ instead of $\Omega$).

\begin{lemma}\label{lemma:thresc1}
\begin{itemize}
\item[(i)] Let $f_0$ and $f^1\in W^{1,2}(\Omega)$ such that $f^h\rightharpoonup f_0$ weakly in $W^{1,2}(\Omega)$ then there exist $\phi\in L^2(\Omega; \dot{W}^{1,2}(\mathcal{Y};\RR^3))$ and $\psi\in L^2(\Omega\times\mathcal{Y};\dot{W}^{1,2}(\mathcal{Y};\RR^{3\times 3}))$ such that
\begin{equation*}
\nabla f^h\overset{3}{\rightharpoonup} \nabla f_0+\nabla_y\phi+\nabla_z\psi\quad\text{weakly $3$-scale}.
\end{equation*}

\item[(ii)] Let $f_0$ and $f^h\in W^{1,2}(\Omega)$ be such that $f^h\rightharpoonup f_0$ weakly in $W^{1,2}(\Omega)$ and assume that $\nabla f^h\overset{3}{\rightharpoonup} \nabla f_0+\nabla_y\phi+\nabla_z\psi$ for some $\phi\in L^2(\Omega; \dot{H}^1(\mathcal{Y};\RR^3))$ and $\psi\in L^2(\Omega\times\mathcal{Y};\dot{W}^{1,2}(\mathcal{Y};\RR^{3\times 3}))$. Then $\displaystyle\frac{f^h}{\varepsilon^2}\overset{osc, Z}{\rightharpoonup} \psi$.

\item[(iii)] Let $f_0$ and $f^h\in W^{1,2}(\Omega)$ be such that $f^h\rightharpoonup f_0$ weakly in $W^{2,2}(\Omega)$ and assume that $\nabla f^h\overset{3}{\rightharpoonup} \nabla f_0+\nabla_y\phi+\nabla_z\psi$ for some $\phi\in L^2(\Omega; \dot{W}^{1,2}(\mathcal{Y};\RR^3))$ and $\psi\in L^2(\Omega\times\mathcal{Y};\dot{W}^{1,2}(\mathcal{Y};\RR^{3\times 3}))$. Moreover, if $\nabla^2 f^h\overset{3}{\rightharpoonup} \nabla^2 f_0+\nabla^2_y\phi+\nabla^2_z\psi$ for some $\phi\in L^2(\Omega; \dot{W}^{2,2}(\mathcal{Y};\RR^3))$ and $\psi\in L^2(\Omega\times\mathcal{Y};\dot{W}^{2,2}(\mathcal{Y};\RR^{3\times 3}))$. Then $\displaystyle\frac{f^h}{\varepsilon^4}\overset{osc, Z}{\rightharpoonup} \psi$.
\end{itemize}
\end{lemma}

\begin{proof}
The proof of (i) can be found in Theorem 1.2 in  \cite{Allamults}. The proof of (ii) is an extension of Lemma 3.7 in \cite{HORVEL14}, and that of (iii) is similar.
\end{proof}

\section{Elasticity framework and intermediate results}\label{sec:3}

Throught this paper we assume that $\varepsilon: (0,1)\to (0,1)$ denotes a function such that the limits
\begin{equation*}
\displaystyle\gamma_1:=\lim_{h\to 0}\frac{h}{\varepsilon(h)}
\end{equation*}
and
\begin{equation*}
\displaystyle\gamma_2:=\lim_{h\to 0}\frac{h}{\varepsilon^2(h)}
\end{equation*}
exist in $[0,+\infty]$. We will frequently write $\varepsilon$ instead of $\varepsilon(h)$, but always with the understanding that $\varepsilon$ depends on $h$. There are five possible regimes: $\gamma_1, \gamma_2=+\infty$; $0<\gamma_1<+\infty$ and $\gamma_2=+\infty$; $\gamma_1=0$ and $\gamma_2=+\infty$; $\gamma_1=0$ and $0<\gamma_2<+\infty$; $\gamma_1=0$ and $\gamma_2=0$. We focus initially on the first three regimes, that is, on the cases in which $\gamma_2=+\infty$, after the case $\gamma_1=0; \gamma_2\in(0,+\infty)$ and then the last one.

\vspace{0.5 cm}

From now on, we consider $\RR^{3\times 3}$: set of all real square matrices of order $3$, the set of all
rotations in $\RR^3$
\begin{equation*}
SO(3)=\left\{\textbf{R}\in \RR^{3\times 3}: \textbf{R}\textbf{R}^T=\textbf{R}^{T}\textbf{R}=\textbf{I},\, \det \textbf{R}=1\right\}
\end{equation*}
and we fix a Borel measurable energy density

\begin{equation*}
W:S^1\times \RR^2\times \RR^2\times \RR^{3\times 3}\to\RR^{+}\cup \{+\infty\}
\end{equation*}
with the following properties:
\begin{itemize}
\item[$\bullet$] $W(\cdot, y, z, \textbf{F})$ is continuous for almost every $y, z\in\RR^2$ and $\textbf{F}\in\RR^{3\times 3}$.
\item[$\bullet$] $W(x,\cdot,\cdot, \textbf{F})$ is $\mathcal{Y}\times\mathcal{Y}-$periodic for all $x\in S^1$ and almost every $\textbf{F}\in\RR^{3\times 3}$.
\item[$\bullet$] For all $(x,y, z)\in S^{1}\times  \mathcal{Y}\times  \mathcal{Y}$ we have $W(x, y, z,I)=0$ and $W(x,y,z, \textbf{R}\textbf{F})=W(x,y,z, \textbf{F})$ for all $\textbf{F}\in\RR^{3\times 3}$, $\textbf{R}\in SO(3)$.
\item[$\bullet$] There exist constants $0<\alpha\leq \beta$ and $\rho>0$ such that for all $(x,y, z)\in S^1\times \mathcal{Y}\times \mathcal{Y}$ we have
\begin{align*}
& \nonumber \hspace{1.5 cm}W(x,y,z,\textbf{F})\geq \alpha\, \text{dist}^2(\textbf{F},SO(3))\quad\forall \textbf{F}\in\RR^{3\times 3}\\
& \nonumber\hspace{1.5 cm} W(x,y,z,\textbf{F})\leq \beta\,\text{dist}^2(\textbf{F},SO(3))\quad\forall \textbf{F}\in\RR^{3\times 3}\,\text{with}\,\text{dist}^2(\textbf{F},SO(3))\leq \rho.
\end{align*}	
\item[$\bullet$] For each $(x,y, z)\in S^1\times \mathcal{Y}\times \mathcal{Y}$ there exists a quadratic form $\mathscr{Q}(x,y, z,\cdot):\RR^{3\times 3}\to \RR$ such that (\ref{eq:linW}) holds.
\end{itemize}
Clearly $\mathscr{Q}(\cdot,y, z,\cdot)$ is continuous for almost every $(y,z)\in\RR^2\times\RR^2$ and $\mathscr{Q}(x,\cdot, \cdot, \textbf{G})$ is $\mathcal{Y}\times\mathcal{Y}-$periodic for all $x\in S^1$ and all $\textbf{G}\in\RR^{3\times 3}$.\\

The elastic energy per unit thickness of a deformation $u^h\in W^{1,2}(S^h;\RR^{3})$ of the shell $S^h$ is given by
\begin{equation*}
J^h(u^h)=\displaystyle\frac{1}{h}\int_{S^h}W\left(\Theta^h(x),r(x)/\varepsilon, r(x)/\varepsilon^2,\nabla u^h(x)\right)dx.
\end{equation*}
In order to express the elastic energy in terms of the new variables, we associate with $y:S^1\to \RR^{3}$ the energy
\begin{align*}
I^h(y)&=\displaystyle\int_{S^1}W(x,r(x)/\varepsilon, r(x)/\varepsilon^2,\nabla_{h}y(x))\det\left(I+t(x)\textbf{S}(x)\right)^{-1}dx\\
&=\displaystyle\int_{S}\int_{I}W(x+tn(x),r(x)/\varepsilon, r(x)/\varepsilon^2,\nabla_{h}y(x+tn(x))\,dt\,d\,vol_{S}(x).
\end{align*}
By a change of variables we have
\begin{equation*}
J^h(u^h)=\displaystyle\frac{1}{h}\int_{S^1}W(x,r(x)/\varepsilon, r(x)/\varepsilon^2, \nabla_{h}y^h(x))\left|\det\nabla(\Theta^h)^{-1}(x)\right|dx,
\end{equation*}
where again $y^h(\Theta^h)=u^h$. Using (\ref{eq:387}) and (\ref{eq:2500087}) we see that there exists a constant $C$ such that
\begin{equation*}
|J^h(u^h)-I^h(y^h)|\leq C h I^h(y^h).
\end{equation*}

\vspace{0.5 cm}

\begin{lemma}\label{lemma:31}
Let $(w^h)\in W^{1,2}(S^1;\RR^3)$ be such that
\begin{equation*} 
\displaystyle\lim_{h\to 0}\sup\left(||w^h||_{L^2(S^1)}+||\nabla_h w^h||_{L^2(S^1)}\right)<\infty.
\end{equation*}
Then there exist a map $w_0\in W^{1,2}(S;\RR^3)$ and a field $H_{\gamma_1, \gamma_2}\in L^2(S\times I\times \mathcal{Y}\times  \mathcal{Y};\RR^{3\times 3})$ of the form
\begin{equation}
H_{\gamma_1, \gamma_2}=
\begin{cases}
(\nabla_{y}w_1+\nabla_{z} w_2, \partial_{3}w_3)& \text{if}\,\, \gamma_1=\gamma_2=+\infty;\partial_{y_i}w_3=\partial_{z_i}w_3\,\,\text{for}\,\, i=1,2,\\
(\nabla_{y}w_1+\nabla_{z} w_2, \partial_{3}w_3)& \text{if}\,\, 0<\gamma_1<+\infty\,\text{and}\,\, \gamma_2=+\infty; w_3=\displaystyle\frac{w_1}{\gamma_1},\\
(\nabla_{y}w_1+\nabla_{z} w_2, \partial_{3}w_3)& \text{if}\,\, \gamma_1=0\,\,\text{and}\,\, \gamma_2=+\infty; \partial_{3}w_1=0\,\,\text{and}\,\, \partial_{z_i}w_3=0,\\
&\,\text{for}\,\, i=1,2,\\
(\nabla_{y}w_1+\nabla_{z} w_2, \partial_{3}w_3)& \text{if}\,\, \gamma_1=\gamma_2=0; \partial_{3}w_1=0\,\,\text{and}\,\, \partial_{3}w_2=0,\\
(\nabla_{y}w_1+\nabla_{z} w_2, \partial_{3}w_3)& \text{if}\,\,\gamma_1=0\,\text{and}\,\, 0<\gamma_2<+\infty; w_3=\displaystyle\frac{w_2}{\gamma_2},
\end{cases}
\end{equation}
for some $w_1\in L^2(S\times I; W^{1,2}(\mathcal{Y}))$, $w_2\in L^2(S\times I\times \mathcal{Y}; W^{1,2}(\mathcal{Y}))$ and $ w_3\in L^2(S\times \mathcal{Y}\times \mathcal{Y};W^{1,2}(I))$, such that, up to a subsequence, $w^h\to w_0$ in $L^2$ and
\begin{equation*}
\nabla_{h}w^h\overset{3}{\rightharpoonup} dw_0\circ T_{S}+\displaystyle\sum_{i,j=1}^{3}(\hat{H}_{\gamma_1, \gamma_2})_{ij}\tau^1\otimes\tau^{j}\quad\text{weakly three-scale on}\, S^1.
\end{equation*}
Here, $\tau^3=n$, $w_0$ is the weak limit in $W^{1,2}(S)$ of $\int_{I}w^h(x+tn(x))dt$ and $\hat{H}_{\gamma_1, \gamma_2}\in L^2(S^1\times  \mathcal{Y}\times  \mathcal{Y};\RR^{3\times 3})$ is defined by $\hat{H}_{\gamma_1, \gamma_2}(x,y, z)=H_{\gamma_1, \gamma_2}(\pi(x), t(x), y, z)$.
\end{lemma}

\vspace{0.5 cm}

\dem
We adapt Lemma 4.3 in \cite{HORVEL15} to our case for $\gamma_1,\gamma_2\in [0,+\infty]$. The hypotheses imply, e.g. by (\ref{eq:2500087}) , that the $w^h$ are uniformly bounded in $W^{1,2}(S^1)$, so up to a subsequence $w^h\rightharpoonup: w_0$ in $W^{1,2}(S^1)$. Set $\tilde{w}^h=w^h\circ \Xi$, so clearly $\tilde{w}^h$ is uniformly bounded in $L^2(\Omega)$. From the uniform $L^2-$bound on $\nabla_{h}w^h$ and from (\ref{eq:2300087}) we deduce that $\tilde{\nabla}_h\tilde{w}^h$ is uniformly bounded in $L^2(\Omega)$. Hence there is $\tilde{w}_0\in W^{1,2}(\Omega;\RR^3)$ with $\partial_{3}\tilde{w}_0=0$ such that $\tilde{w}^h\rightharpoonup \tilde{w}_0$ weakly in $W^{1,2}(\Omega;\RR^3)$; clearly $\tilde{w}_0=w_0\circ\Xi$, so (since $\partial_{3}\tilde{w}_0$=0) in particular $w_0$ is the trivial extension of a map defined on $S$.
 \medskip In the case $\gamma_{1}=0$ and $\gamma_2=+\infty$ by uniform boundedness in $L^2(\Omega)$, there exist (see Theorem 3.2 in \cite{BUFDAVFON15}) $\tilde{w}_1\in L^2(\Omega;W^{1,2}(\mathcal{Y})))$, $\tilde{w}_2\in L^2(\Omega\times \mathcal{Y}; W^{1,2}(\mathcal{Y}))$ and $\overline{w}_3\in L^2(\omega\times \mathcal{Y}\times \mathcal{Y}; W^{1,2}(I))$ with $\partial_{y_i}\overline{w_3}=\partial_{z_i}\overline{w_3}=0$ for $i=1,2$ such that, up to the extraction of a (not relabeled) subsequence,
 \begin{equation*}
 \tilde{\nabla}_h\tilde{w}^h\overset{3}{\rightharpoonup} (\partial_{1}\tilde{w}_0, \partial_{2}\tilde{w}_0, 0)+\left(\partial_{y_1}\tilde{w}_1, \partial_{y_2}\tilde{w}_1, 0\right)+\left(\partial_{z_1}\tilde{w}_2, \partial_{z_2}\tilde{w}_2, \partial_{3}\overline{w}_3\right)\quad \text{in}\,\,\Omega.
 \end{equation*}
By (\ref{eq:2300087}) the left-hand side equals $(\nabla_h w)(\Xi)\nabla\Xi(\tilde{\Theta}^{-1}_{h})$. As $\nabla\Xi(\tilde{\Theta}^{-1}_{h})$ converges uniformly on $S^1$ to $(\partial_1\xi, \partial_2\xi, n(\xi))$ (extended trivially in the $x_3-$direction), we conclude:
\begin{equation*}
 \nabla_hw^h(\Xi)\overset{3}{\rightharpoonup} \left((\partial_{1}\tilde{w}_0, \partial_{2}\tilde{w}_0, 0)+\left(\partial_{y_1}\tilde{w}_1, \partial_{y_2}\tilde{w}_1, 0\right)+\left(\partial_{z_1}\tilde{w}_2, \partial_{z_2}\tilde{w}_2, \partial_{3}\overline{w}_3\right)\right)(\partial_1\xi, \partial_2\xi, n(\xi))^{-1}.
  \end{equation*}
  On the right-hand side we use
  \begin{equation*}
  (\partial_1\xi, \partial_2\xi, n(\xi))^{-1}\circ \Xi^{-1}=(\tau_1, \tau_2, n)^{-1}=(\tau^1, \tau^2, n)^{T}
  \end{equation*}
  and $(\partial_{\alpha}\tilde{w}_0)\circ \Xi^{-1}=d w_0(\tau_{\alpha})$ to obtain the claim when  $\gamma_{1}=\gamma_2=+\infty$, after defining $(w_1)_{i}=(\tilde{w}_1\circ r)\cdot \tau_{i}$, $(w_2)_{i}=(\tilde{w}_2\circ r)\cdot \tau_{i}$  and $(w_3)_{i}=(\overline{w}_3\circ r)\cdot \tau_{i}$ for $i=1,2, 3$. The other four cases are proven similarly using Theorem 3.2 in \cite{BUFDAVFON15} in where we can cover the two last cases ($ \gamma_1=\gamma_2=0$; $\gamma_1=0\,\text{and}\,\, 0<\gamma_2<+\infty$) which are valid for plates and can naturally be extended to the case of shells.
\lqqd

\subsection{Asymptotic energy functionals}\label{subsectiongammas}

Next we will introduce the asymptotic energy functionals. To do so, we need the definition of the relaxation fields and the cell formulae. Recall that $a\odot b=\frac{1}{2}(a\otimes b+b\otimes a)$. We make the following definitions:

\begin{equation*}
D(\mathcal{U}_{0, +\infty})=W^{1,2}(\mathcal{Y},\RR^2)\times L^2(I\times\mathcal{Y}; W^{1,2}(\mathcal{Y}; \RR^3))\times W^{2,2}(\mathcal{Y})\times L^2(I\times\mathcal{Y},\RR^3)
\end{equation*}

 and for $(\zeta, \eta, \varphi, \mu)\in L^2(S, D(\mathcal{U}_{0, +\infty}))$ define
\begin{equation*}
\mathcal{U}_{0, +\infty}(\zeta,\eta, \varphi, \mu)=\text{Def}_{\mathcal{Y}}\,\zeta+2\mu_{\alpha}\tau^{\alpha}\odot n+\mu_{3} n\odot n-t\text{Hess}_{\mathcal{Y}}\,\varphi+\text{Def}_{\mathcal{Z}}\,\eta+2\partial_{z_{\alpha}}\eta_3\tau^{\alpha}\odot n.
\end{equation*}

\begin{equation*}
D(\mathcal{U}_{+\infty, +\infty})=L^2(I, W^{1,2}(\mathcal{Y},\RR^2))\times L^2(I\times\mathcal{Y}; W^{1,2}(\mathcal{Y}; \RR^3))\times L^2(I, W^{1,2}(\mathcal{Y}))\times L^2(I,\RR^3)
\end{equation*}
and for $(\zeta, \eta, \rho, c)\in L^2(S, D(\mathcal{U}_{+\infty, +\infty})$ define

\begin{equation*}
\mathcal{U}_{+\infty, +\infty}(\zeta,\eta, \rho, c)=\text{Def}_{\mathcal{Y}}\,\zeta+2(\partial_{y_{\alpha}}\rho+c_{\alpha})\tau^{\alpha}\odot n+c_3n\odot n+\text{Def}_{\mathcal{Z}}\,\eta++2\partial_{z_{\alpha}}\eta_3\tau^{\alpha}\odot n.
\end{equation*}
For $\gamma_{1}\in(0,+\infty)$ 
\begin{equation*}
D(\mathcal{U}_{\gamma_1, +\infty})=W^{1,2}(I\times\mathcal{Y},\RR^2)\times L^2(I\times \mathcal{Y}; W^{1,2}(\mathcal{Y};\RR^3))\times W^{1,2}(I\times\mathcal{Y})
\end{equation*}
and for $(\zeta, \eta,\rho)\in L^2(S, D(\mathcal{U}_{\gamma_1, +\infty}))$ define
\begin{equation*}
\mathcal{U}_{\gamma_1, +\infty}(\zeta,\eta, \rho)=\displaystyle\text{Def}_{\mathcal{Y}}\,\zeta+\left(\partial_{y_{\alpha}}\rho+\frac{1}{\gamma_1}\partial_{3}\zeta_{\alpha}\right)\tau^{\alpha}\odot n+\left(\frac{1}{\gamma_1}\partial_{3}\rho\right)n\odot n+\text{Def}_{\mathcal{Z}}\,\eta++2\partial_{z_{\alpha}}\eta_3\tau^{\alpha}\odot n.
\end{equation*}

By embedding $D(\mathcal{U}_{0,+\infty})$ trivially into $L^2(S, D(\mathcal{U}_{0,\infty}))$, we can regard $\mathcal{U}_{0,+\infty}$ as a map from $D(\mathcal{U}_{0,+\infty})$ into $L^2(S, L^2(I\times\mathcal{Y}\times\mathcal{Y}, \RR^{3\times 3}_{\text{sym}}))$.

For $\gamma_2\in (0,+\infty)$ set 
\begin{equation*}
D(\mathcal{U}_{0, \gamma_2})=W^{1,2}(\mathcal{Y},\RR^2)\times L^2(I\times\mathcal{Y}; W^{1,2}(\mathcal{Y}; \RR^3))\times W^{2,2}(\mathcal{Y})\times L^2(I\times\mathcal{Y},\RR^3)
\end{equation*}

 and for $(\zeta, \eta, \varphi, \mu)\in L^2(S, D(\mathcal{U}_{0,\gamma_2}))$ define
\begin{equation*}
\mathcal{U}_{0, \gamma_2}(\zeta,\eta, \varphi, \mu)=\text{Def}_{\mathcal{Y}}\,\zeta+2\mu_{\alpha}\tau^{\alpha}\odot n+\mu_{3} n\odot n-t\text{Hess}_{\mathcal{Y}}\,\varphi+\displaystyle\frac{1}{\gamma_2}\varphi\textbf{S}(x)+\text{Def}_{\mathcal{Z}}\,\eta+2\partial_{z_{\alpha}}\eta_3\tau^{\alpha}\odot n.
\end{equation*}

For each $x\in S$ the fiberwise action $\mathcal{U}^{(x)}_{0, +\infty}$ of $\mathcal{U}_{0, +\infty}$ is
\begin{align*}
\mathcal{U}^{(x)}_{0, +\infty}(\zeta, \eta, \varphi, \mu)&=(\text{Def}_{\mathcal{Y}}\,\zeta)(x)+2\mu_{\alpha}\tau^{\alpha}(x)\odot n(x)+\mu_{3} n(x)\odot n(x)-t(\text{Hess}_{\mathcal{Y}}\,\varphi)(x)\\
&\quad+(\text{Def}_{\mathcal{Z}}\,\eta)(x)+2\partial_{z_{\alpha}}\eta_3\tau^{\alpha}\odot n,
\end{align*}
for all $(\zeta, \eta, \varphi, \mu)\in D(\mathcal{U}_{0, +\infty})$.\\

\medskip For each $x\in S$ we define $L^{(x)}_{0, +\infty}(I\times \mathcal{Y}\times\mathcal{Y})=\mathcal{U}^{(x)}_{0, +\infty}(D(\mathcal{U}_{0, +\infty}))$, i.e.,
\begin{equation*}
L^{(x)}_{0, +\infty}(I\times\mathcal{Y}\times\mathcal{Y})=\left\{\mathcal{U}^{(x)}_{0, +\infty}(\zeta, \eta, \varphi, \mu) : (\zeta, \eta, \varphi, \mu)\in D(\mathcal{U}_{0, +\infty})\right\}.
\end{equation*}
This is a subspace of $L^2(I\times\mathcal{Y}\times\mathcal{Y}, \RR^{3\times 3}_{\text{sym}})$. We denote by $L_{0, +\infty}(I\times\mathcal{Y}\times\mathcal{Y})$ the vector bundle over $S$ with fibers $L^{(x)}_{0, +\infty}(I\times\mathcal{Y}\times\mathcal{Y})$; in what follows we will frequently omit the index $(x)$ for the fibers. The bundle $L_{\gamma_1, +\infty}(I\times \mathcal{Y}\times\mathcal{Y})$, for $\gamma_{1}\in (0,+\infty]$ and $L_{0,\gamma_2}(I\times \mathcal{Y}\times\mathcal{Y})$, for $\gamma_{2}\in (0,+\infty)$ are defined analogously.  The elements of these spaces are the relaxation fields.

\medskip For $\gamma_1\in[0,+\infty]$, $\gamma_2=+\infty$ and $x\in S$, we define $\mathscr{Q}_{\gamma_1, +\infty}(x,\cdot);T^{*}_{x}S\otimes T^{*}_{x}S\to \RR$ by setting
\begin{equation*}
\mathscr{Q}_{\gamma_1, +\infty}(x, q)=\inf \int_{I}\int_{\mathcal{Y}}\int_{\mathcal{Y}}\mathscr{Q}\left(x+tn(x), y, z, p+tq+U(t,y,z)\right)\,dz\,dy\,dt.
\end{equation*}
Here the infimum is taken over all $U\in L^{(x)}_{\gamma_1, +\infty}(I\times\mathcal{Y}\times\mathcal{Y})$ and all $p\in T^{*}_{x}S\otimes T^{*}_{x}S$.

\medskip For $\gamma_1=0$, $\gamma_2\in(0,+\infty)$ and $x\in S$, we define $\mathscr{Q}_{0, \gamma_2}(x,\cdot);T^{*}_{x}S\otimes T^{*}_{x}S\to \RR$ by setting
\begin{equation*}
\mathscr{Q}_{0,\gamma_2}(x, q)=\inf \int_{I}\int_{\mathcal{Y}}\int_{\mathcal{Y}}\mathscr{Q}\left(x+tn(x), y, z, p+tq+U(t,y,z)\right)\,dz\,dy\,dt.
\end{equation*}
Here the infimum is taken over all $U\in L^{(x)}_{0, \gamma_2}(I\times\mathcal{Y}\times\mathcal{Y})$ and all $p\in T^{*}_{x}S\otimes T^{*}_{x}S$.

\medskip Note that $\mathscr{Q}_{\gamma_{1}, +\infty}(x, q)=\mathscr{Q}_{\gamma_1, +\infty}(x, \text{sym}\, q)$ for all $x\in S$ and all $q\in T^{*}_{x}S\otimes T^{*}_{x}S$. For $x\in S$ and $q\in T^{*}_{x}S\odot T^{*}_{x}S$ define the homogeneous relaxation (cf. \cite{Lewicakamorapa}):
\begin{equation*}
\widetilde{\mathscr{Q}}(x, t, q)=\min_{M\in\RR^{3\times 3}_{\text{sym}}}\left\{\mathscr{Q}(x+tn(x), M): M(T_{S}, T_{S})=q(T_S, T_{S})\right\}.
\end{equation*}
Then it is easy to see that
\begin{align*}
&\mathscr{Q}_{0, +\infty}(x, q)=\inf \displaystyle\int_{I\times\mathcal{Y}\times\mathcal{Y}}\widetilde{\mathscr{Q}}\left(x, y, z, p+tq+(\text{Def}_{\mathcal{Y}}\zeta)(x)-t(\text{Hess}_{\mathcal{Y}}\varphi)(x)+(\text{Def}_{\mathcal{Z}}\eta)(x)\right)\,dz\, dy\, dt,\\
&\mathscr{Q}_{0, \gamma_2}(x, q)\\
&=\inf \displaystyle\int_{I\times\mathcal{Y}\times\mathcal{Y}}\widetilde{\mathscr{Q}}\left(x, y, z, p+tq+(\text{Def}_{\mathcal{Y}}\zeta)(x)-t(\text{Hess}_{\mathcal{Y}}\varphi)(x)+\gamma_2^{-1}\varphi\,\textbf{S}(x)+(\text{Def}_{\mathcal{Z}}\eta)(x)\right)\,dz\, dy\, dt,
\end{align*}

where the infimum is taken over all $\zeta\in W^{1,2}(\mathcal{Y},\RR^2)$, all $\varphi\in W^{2,2}(\mathcal{Y})$, all $\eta\in L^2(I\times\mathcal{Y}; W^{1,2}(\mathcal{Y}; \RR^3))$ and all $p\in T^{*}_{x}S\odot T^{*}_{x}S$. In the case when the material is homogeneous in the thickness direction, we have
\begin{equation*}
\mathscr{Q}_{0, +\infty}(x, q)=\displaystyle\frac{1}{12}\inf\left\{\int_{\mathcal{Y}}\int_{\mathcal{Y}}\widetilde{Q}(x,y, z, q+(\text{Hess}_{\mathcal{Y}}\,\varphi)(x))\,dz\,dy : \varphi\in \dot{W}^{2,2}(\mathcal{Y})\right\}.
\end{equation*}

The analogous formula holds for $\mathscr{Q}_{0, \gamma_2}$.

\medskip As in \cite{HORVEL18}, for all $x\in S$ and all $q\in T^{*}_{x}S\odot T^{*}_{x}S$ we have
\begin{equation*}
\displaystyle\lim_{\gamma_{1}\to \infty}\mathscr{Q}_{\gamma_1, +\infty}(x, q)=\mathscr{Q}_{+\infty, +\infty}(x, q)\quad\text{and}\quad\lim_{\gamma_1\to 0}\mathscr{Q}_{\gamma_1, +\infty}(x, q)=\mathscr{Q}_{0, +\infty}(x, q).
\end{equation*}
It is not difficult to show that for all $\gamma_1\in[0, +\infty]$, $\gamma_2=+\infty$ and $x\in S$ the map $\mathscr{Q}_{\gamma_1, +\infty}(x, \cdot)$ is quadratic and that there exist $c_1, c_2>0$ such that for all $x\in S$ we have
\begin{equation*}
c_1 |\text{sym}\, q|^2\leq \mathscr{Q}_{\gamma_1, +\infty}(x, q)\leq c_2|\text{sym}\, q|^2,\quad \forall q\in T^{*}_{x}S\otimes T^{*}_{x}S.
\end{equation*}

For $\gamma_1\in[0,+\infty]$ and $\gamma_2=+\infty$ we define $I_{\gamma_1, +\infty}:W^{1,2}(S;\RR^3)\to\RR$ by setting
\begin{equation*}
I_{\gamma_1, +\infty}(u)=\begin{cases}
\int_{S}\mathscr{Q}_{\gamma_1, +\infty}(x, \textbf{S}^{r}_{u}(x))\,d\,\text{vol}_{S}(x) &\text{if}\, u\in W^{2,2}_{\text{iso}}(S), \\
+\infty &\text{otherwise}. 
\end{cases}
\end{equation*}

For $\gamma_2\in(0,+\infty)$ the functional $I_{0, \gamma_2}:W^{1,2}(S;\RR^3)\to\RR$ is defined analogously, by replacing $\mathscr{Q}_{\gamma_1, +\infty}$ by $\mathscr{Q}_{0, \gamma_2}$.

\section{Main result}\label{sec:4}

For a given sequence $(u^h)\subset W^{1,2}(S^h;\RR^3)$ we continue to define the sequence $(y^h)\subset W^{1,2}(S^1,\RR^3)$ of rescaled deformations by $y^h(\Theta^h)=u^h$. We recall the compactness result for sequences with finite bending energy, cf. Theorem 1 in \cite{FJM06} for a proof.

\vspace{0.5 cm}

\begin{prop}\label{prop:41}
Let $(u^h)\subset W^{1,2}(S^h,\RR^3)$ be such that (\ref{eq:33}) holds. Then there exists $u\in W^{2,2}_{\text{iso}}(S)$ such that (after passing to subsequences and extending $u$ and $n$ trivially to $S^1$) as $h\to 0$ we have
\begin{align}
y^h-\displaystyle\frac{1}{|S^1|}\int_{S^1}y^hdx&\label{may820231}\to u\,\,\text{strongly in}\,\, W^{1,2}(S^1;\RR^3),\\
\nabla_{h}y^h&\label{may820232}\to Q\,\,\text{strongly in }\,\, L^2(S^1;\RR^{3\times 3}).
\end{align}
Here $Q\in W^{1,2}(S,\text{SO}(3))$ is determined by the condition $Q\tau=\nabla_{\tau}u$ for all smooth tangent vector fields $\tau$ along $S$.
\end{prop}

We denote by $\widetilde{W}^{2,2}_{\text{iso}} (S)$
 the set of those maps $u\in W^{2,2}_{\text{iso}}(S)$ 
  for which there exists $(u^h)\subset W^{2,\infty}_{\text{iso}}(S)$ converging strongly to $u$ in $W^{2,2}(S)$. The reason for introducing this space is that we can construct the recovery sequence only to limit deformations $u$ belonging to this space. Theorem 2.1 in \cite{HORVEL18} plays an essential role in this construction. The following $\Gamma-$convergence result is our main result:
\begin{teo}\label{teo:187}
Let $\gamma_{1}\in[0,+\infty]$ and $\gamma_2=+\infty$. Then the following are true:
\begin{itemize}
\item[(i)] Let $(u^h)\subset W^{1,2}(S^h,\RR^3)$ be such that (\ref{eq:33}) holds and such that $y^h-\frac{1}{|S^1|}\int_{S^1}y^h\to u$ strongly in $L^2(S^1)$ for some $u\in L^2(S^1, \RR^3)$. Then
\begin{equation*}
\lim_{h\to 0}\inf h^{-2} J^h(u^h)\geq I_{\gamma_1, +\infty}(u).
\end{equation*}
\item[(ii)] If, in addition, $S$ is simply connected, then for every $u\in \widetilde{W}^{2, 2}_{\text{iso}}(S)$ there exists $(u^h)\subset W^{1,2}(S^h;\RR^3)$ satisfying (\ref{eq:33}), and such that $y^h\to u$ strongly in $W^{1,2}(S^1)$. Moreover,
\begin{equation*}
\lim_{h\to 0}h^{-2}J^h(u^h)=I_{\gamma_1,+\infty}(u).
\end{equation*}
\end{itemize}
Furthermore, for $\gamma_1=0$ and $\gamma_2\in(0,+\infty)$ the items $(i)$ and $(ii)$ hold replacing $I_{\gamma_1, +\infty}$ by $I_{0, \gamma_2}$.
\end{teo}

\section{Proof of main result}\label{sec:5}

We present the proof of Theorem \ref{teo:187} in the following two subsections: 
\subsection{Proof of Theorem \ref{teo:187} i). }\label{subsect51}
The next results are important tools for proving Theorem \ref{teo:187} i). 

 We consider a sequence $(u^h)\subset W^{1,2}(S^h, \RR^3)$ satisfying (\ref{eq:33}) and we set $y^h(\Theta^h)=u^h$.

\begin{lemma}\label{lemma:33}
Define
\begin{equation*}
\delta=
\begin{cases}
\varepsilon, & \text{if}\,\, \gamma_1\in(0,+\infty), \gamma_2=+\infty,\\
\displaystyle\left(2\lceil\frac{h}{\varepsilon}\rceil+1\right)\varepsilon, &\text{if}\,\,\gamma_1=\gamma_2=+\infty,\\
\displaystyle\left(2\lceil\frac{h}{\varepsilon^2}\rceil+1\right)\varepsilon^2,  &\text{if}\,\, \gamma_1=0, \gamma_2=+\infty,\\
\varepsilon^2, & \text{if}\,\,\gamma_1=0, \gamma_2\in(0,+\infty),\\
h, & \text{if}\,\, \gamma_1=\gamma_2=0.
\end{cases}
\end{equation*}
Then there exist constants $C, c>0$ such that the following is true: if $h\leq c$ and $u\in W^{1,2}(S^h;\RR^3)$, then there exists a map $\tilde{R}:\omega\to SO(3)$ which is constant on each cube $x+\delta Y$ with $x\in\delta \ZZ$ and there exist $\tilde{R}_{s}\in W^{1,2}(\omega;\RR^3)$ such that for each $a\in\RR^2$ with $|a_1|\leq \delta$ and $|a_2|\leq \delta$ and for each $\tilde{\omega}\subset\omega$ with $\text{dist}\,(\tilde{\omega}, \partial\omega)>c\delta$ we have:
\begin{align}
&\nonumber ||(\nabla_{h}y)(\Xi)-\tilde{R}||^2_{L^2(\tilde{\omega}\times I)}+||\tilde{R}-\tilde{R}_{s}||^2_{L^2(\tilde{\omega})}+h^2||\tilde{R}-\tilde{R}_{s}||^2_{L^{\infty}(\tilde{\omega})}\\
&\nonumber+h^2||(\partial_{1}\tilde{R}_{s}, \partial_{2}\tilde{R}_{s})||^2_{L^2(\tilde{\omega})}+||\tilde{R}(\cdot+a)-\tilde{R}||^2_{L^2(\tilde{\omega})}\\
&\label{constructrotation}\leq C\displaystyle\int_{\Omega}\text{dist}^2(\nabla_{h}y(\Xi), SO(3)).
\end{align}
\end{lemma}

\begin{proof}
This lemma is essentially contained in \cite{BUFDAVFON15}; it is direct consequence of Theorem 3.1 in \cite{FJM02} and arguments in \cite{HORVEL18}, so we will limit ourselves to provide an outline and brief indications on how to construct it. For the first three cases $\gamma_1\in [0,+\infty]$ and $\gamma_2=+\infty$ (i.e., $\gamma_1\in(0,+\infty), \gamma_2=+\infty$, $ \gamma_1=\gamma_2=+\infty$ and $\gamma_1=0, \gamma_2=+\infty$)  we apply the lemma \ref{lemma:a2} with $\delta(h)=\varepsilon(h), \displaystyle\left(2\left\lceil\frac{h}{\varepsilon(h)}\right\rceil+1\right)\varepsilon(h)$ and $\displaystyle\left(2\left\lceil\frac{h}{\varepsilon^2(h)}\right\rceil+1\right)\varepsilon^2(h)$, respectively.  The novelty in this work lies in the cases $\gamma_1=\gamma_2=0$ and $\gamma_1=0, \gamma_2\in(0,+\infty)$, when similarly we apply the lemma \ref{lemma:a2} with $\delta(h)=h$ and $\delta(h)=\varepsilon^2(h)$, respectively. Finally, we can rewrite this previously obtained result as an analogous of the Lemma \ref{lemma:a2}. The motivation for choosing these $\delta$ explicitly follows from the way sequences $\{R^h\}\subset L^{\infty}(\omega; SO(3))$ and $\{\tilde{R}^h\}\subset W^{1,2}(\omega; \RR^{3\times 3})$ such that is piecewise constant on every cube of the form $Q(\varepsilon(h)z,\varepsilon(h))$ with $z\in \ZZ^2$ such that (\ref{rotationRandRtilde}) holds, will be constructed, aiming to identify the three-scale limit of scaled linearized stresses (\ref{eq:34}). Such choice of such $\delta$ in each case will be justified in the proof of the following proposition. 

\end{proof}


\begin{prop}\label{prop:5}
Let $\gamma_1\in(0,+\infty)$ and $\gamma_2=+\infty$. Let $(u^h)\subset W^{1,2}(S^h;\RR^3)$ satisfying (\ref{eq:33}) and let $u\in W^{2,2}_{\text{iso}}(S)$ such that (\ref{may820231}) and (\ref{may820232}) hold. Let $\tilde{\omega}\subset\RR^2$ be a domain with $C^{1,1}$ boundary whose closure is contained in $\omega$ and set $\widetilde{S}=\xi(\tilde{\omega})$.
Denote by $\tilde{R}^h:\omega\to SO(3)$ the piecewise constant map obtained by applying Lemma \ref{lemma:33} to $u^h$ and define $R^h:S^1\to SO(3)$ by $R^h=\tilde{R}^h\circ r$. Define $G^h\in L^2(S^1;\RR^{3\times 3})$ by (\ref{eq:34}), where $y^h(\Theta^h)=u^h$. Then there exist $B\in L^2(\tilde{S}, T^{*}\tilde{S}\odot T^{*}\tilde{S})$ and $(\zeta, \eta, \rho)\in L^2(\tilde{S}, D(\mathcal{U}_{\gamma_1, +\infty}))$ such that (up to a subsequence) 
\begin{equation}\label{eq:35}
\text{sym}\,G^h\overset{3}{\rightharpoonup} B+t\textbf{S}^{r}_{u}+\mathcal{U}_{\gamma_1, +\infty}(\zeta, \eta, \rho).
\end{equation}
If $\gamma_1=+\infty$ and $\gamma_2=+\infty$, (\ref{eq:35}) holds true with $\mathcal{U}_{+\infty, +\infty}(\zeta, \eta, \rho, c)$ in place of  $\mathcal{U}_{\gamma_1, +\infty}(\zeta, \eta, \rho)$, where $(\zeta, \eta, \rho, c)\in L^2(\widetilde{S},D(\mathcal{U}_{+\infty, +\infty}))$. Moreover, if $\gamma_1=0$ and $\gamma_2=+\infty$, (\ref{eq:35}) holds true  with $\mathcal{U}_{0, +\infty}(\zeta, \eta, \varphi, \mu)$ in place of $\mathcal{U}_{\gamma_1, +\infty}(\zeta, \eta, \rho)$, where $(\zeta, \eta, \varphi, \mu)\in L^2(\widetilde{S},D(\mathcal{U}_{0, +\infty}))$. Likewise for $\gamma_1=0$ and $\gamma_2\in(0,+\infty)$, (\ref{eq:35}) holds true  with $\mathcal{U}_{0, \gamma_2}(\zeta, \eta, \varphi, \mu)$ in place of $\mathcal{U}_{\gamma_1, +\infty}(\zeta, \eta, \rho)$, where $(\zeta, \eta, \varphi, \mu)\in L^2(\widetilde{S},D(\mathcal{U}_{0, \gamma_2}))$. We have convergence (\ref{eq:35}), also, for $\gamma_1=\gamma_2=0$ but this time, we get $\mathcal{U}_{0, 0}(\dot{B}, \eta, \mu)$ in place of $\mathcal{U}_{\gamma_1, +\infty}(\zeta, \eta, \rho)$, where $(\dot{B}, \eta, \mu)\in L^2(\widetilde{S},D(\mathcal{U}_{0, 0}))$.
\end{prop}
\begin{proof}
Define $\overline{u}^h: S\to\RR^3$ by setting
\begin{equation*}
\overline{u}^h(x)=\displaystyle\frac{1}{h}\int_{hI}u^h(x+tn(x))dt\quad\text{for all}\,\,x\in S.
\end{equation*}

Let $\tilde{R}^h_s:\tilde{\omega}\to\RR^{3\times 3}$ be the maps obtained by applying Lemma \ref{lemma:33} to $u^h$ and set $R^h_{s}=\tilde{R}^h_{s}\circ r$ and $R^h=\tilde{R}^h\circ r$. On $\tilde{S}^h$ define $z^h$ via
\begin{equation*}
u^h=\overline{u}^h(\pi)+t(R^h_{s}n)(\pi)+hz^h.
\end{equation*}
Clearly
\begin{equation}\label{normaleq:april18}
\nabla_{n}u^h=(R^h_{s}n)(\pi)+h\nabla_{n}z^h.
\end{equation}
Let $\tau$ be a smooth tangent vector field along $S$. Then we have
\begin{equation}\label{tangenteq:april18}
\nabla_{\tau}u^h=\nabla_{\nabla_{\tau}\pi}\overline{u}^h(\pi)+t(\nabla_{\nabla_{\tau}\pi}R^h_{s})(\pi)n(\pi)+t(R^h_{s}\textbf{S})(\pi)\nabla_{\tau}\pi+h\nabla_{\tau}z^h.
\end{equation}
Observe that (\ref{eq:281}) implies that $\nabla_{\tau}\pi$ equals $\tau-t\textbf{S}\tau$ up to a term of higher order. Using this and rewriting the problem in coordinates, one can now argue as in Theorem 4.1 in \cite{BUFDAVFON15} to deduce the claim for $\gamma_1\in[0,+\infty]$. We will outline the proof in each case as follows:\\

\textbf{Case 1:} $\gamma_1\in (0,+\infty)$ and $\gamma_2=+\infty$.\\

In order to compute the three-scale limit of
\begin{equation*}
G^h=\displaystyle\frac{(R^h)^T\nabla_h y^h-I}{h},
\end{equation*}
we can rewrite as follow
\begin{align}
R^hG^h&\nonumber=\displaystyle\frac{\nabla_h y^h-R^h}{h}=\frac{\nabla_{\tau}u^h-(R^h\tau)(\pi)}{h}+\frac{h^{-1}\nabla_{n}u^h-(R^hn)(\pi)}{h},
\end{align}
therefore by virtue of (\ref{normaleq:april18})and (\ref{normaleq:april18}) we obtain
\begin{align}
R^hG^h&\nonumber=\displaystyle\frac{\nabla_{\nabla_{\tau}\pi}\overline{u}^h(\pi)-(R^h\tau)(\pi)}{h}+\frac{t}{h}(\nabla_{\nabla_{\tau}\pi}R^h_{s})(\pi)n(\pi)+\frac{t}{h}(R^h_{s}\textbf{S})(\pi)\nabla_{\tau}\pi\\
&\label{splitGh}\quad+\displaystyle\frac{(R^h_{s}n)(\pi)-(R^hn)(\pi)}{h}+\nabla_{h}z^h.
\end{align}
The main difficulty lies in providing a characterization of the three-scale limit for the first term of (\ref{splitGh}) which we can call V, i.e.
\begin{equation}\label{3scaleV}
\displaystyle\frac{\nabla_{\nabla_{\tau}\pi}\overline{u}^h(\pi)-(R^h\tau)(\pi)}{h}\overset{3}{\rightharpoonup} V\quad \text{weakly three-scale in}\,\,L^2(S^1).
\end{equation}
By application of remark \ref{remarkthreescaleconvmarch21} then (\ref{3scaleV}) is equivalent to
\begin{equation}\label{3scaleV2}
\displaystyle\frac{\nabla'\overline{y}^h-(\tilde{R}^h)'}{h}\overset{3}{\rightharpoonup} \tilde{V}\quad \text{weakly three-scale in}\,\,L^2(\Omega).
\end{equation}
Following closely the proof of Theorem 4.1  in \cite{BUFDAVFON15}, this provides us with the characterization now of $\tilde{V}$, which is conditioned on proving over the plate $\Omega$ that 
\begin{equation}\label{conditiontoprovemach21}
\displaystyle\int_{\omega}\frac{(\tilde{R}^h(x')}{h}:(\nabla')^{\perp}\varphi^{\varepsilon}(x')\psi\left(x',\frac{x'}{\varepsilon(h)}\right)dx'=0,
\end{equation}
where $\varphi^{\varepsilon}=\displaystyle\varepsilon^2(h)\varphi\left(\frac{x'}{\varepsilon^2(h)}\right)$.
In order to prove (\ref{conditiontoprovemach21}) the Friesecke, James and M\"uller's rigidity estimate (Theorem 4.1 \cite{FJM02}) provides sequences of rotations $\{R^h\}$ that are piecewise constant on cubes of sizes $\varepsilon(h)$ with center in $\varepsilon(h)\ZZ^2$, however the sequence of test functions $\varphi^{\varepsilon}$ in (\ref{conditiontoprovemach21}) have oscillations on scale $\varepsilon^2(h)$ . This previous problem related to the size oscillations is solved as in theorem 4.2.1 in \cite{BUFDAVFON15}, by subdividing the cubes of sizes $\varepsilon^2(h)$, with centers in $\varepsilon^2(h)\ZZ^2$, into ``good cubes" lying completely within a bigger cube of size $\varepsilon(h)$ and center $\varepsilon(h)\ZZ^2$ and ``bad cubes", and by showing that the measure of the intersection between $\omega$ and the set of ``bad cubes" converges to zero faster with than are comparable to $\varepsilon(h)$, as $h\to 0$ (for details see Remark \ref{remarkcase1}). For the three-scale convergence of the remaining terms in (\ref{splitGh}), we apply the Theorem 1.2 in \cite{Allamults} and Remark \ref{oscilationandthrescale}.The fields $\mathcal{U}_{\gamma_1, +\infty}$ arise, essentially, due the Lemma \ref{lemma:31}.\\

\textbf{Case 2:} $\gamma_1=+\infty$ and $\gamma_2=+\infty$.\\

The proof is similar to the previous case, therefore we only outline the main modifications. Let $\tilde{R}^h_s:\tilde{\omega}\to\RR^{3\times 3}$ be the maps obtained by applying Lemma \ref{lemma:33} with $\delta(h)=\displaystyle\left(2\left\lceil\frac{h}{\varepsilon(h)}\right\rceil+1\right)\varepsilon(h)$ to $u^h$ and set $R^h_{s}=\tilde{R}^h_{s}\circ r$ and $R^h=\tilde{R}^h\circ r$, where $\lceil s\rceil$ denotes the smallest positive integer larger or equal to $s$. 
By the choice of $\delta$ we have
\begin{equation*}
\displaystyle\lim_{h\to 0}\frac{h}{\delta(h)}=\frac{1}{2}
\end{equation*}
and the maps $R^h$ are piecewise constant on cubes of the form $Q(\delta(h)z, \delta(h))$, with $z\in\ZZ^2$. Since $\left\{\frac{\delta(h)}{\varepsilon(h)}\right\}$ is a sequence of odd integers, by virtue of Lemma \ref{Lemma:a3} the maps $R^h$ are piecewise constant on cubes of the form $Q(\varepsilon(h) z, \varepsilon (h))$ with $z\in\ZZ$ and the estimates in Lemma \ref{lemma:33} holds true. Finally, arguing as in the previous case (considering ``good cubes" and ``bad cubes") we obtain the characterization of $G^h$. The fields $\mathcal{U}_{+\infty, +\infty}$ arise, essentially, due the Lemma \ref{lemma:31}.\\

\textbf{Case 3:} $\gamma_1=0$ and $\gamma_2=+\infty$.\\

The proof is similar to the previous case, therefore we only outline the main modifications. Let $\tilde{R}^h_s:\tilde{\omega}\to\RR^{3\times 3}$ be the maps obtained by applying Lemma \ref{lemma:33} with $\delta(h)=\displaystyle\left(2\left\lceil\frac{h}{\varepsilon^2(h)}\right\rceil+1\right)\varepsilon^2(h)$ to $u^h$ and set $R^h_{s}=\tilde{R}^h_{s}\circ r$ and $R^h=\tilde{R}^h\circ r$. Again by the choice of $\delta$ we have
\begin{equation*}
\displaystyle\lim_{h\to 0}\frac{h}{\delta(h)}=\frac{1}{2}
\end{equation*}
and the maps $R^h$ are piecewise constant on cubes of the form $Q(\delta(h)z, \delta(h))$, with $z\in\ZZ^2$. Since $\left\{\frac{\delta(h)}{\varepsilon^2(h)}\right\}$ is a sequence of odd integers, by virtue of Remark \ref{remarkA2} the maps $R^h$ are piecewise constant on cubes of the form $Q(\varepsilon^2(h) z, \varepsilon^2 (h))$ with $z\in\ZZ$ and the estimates in Lemma \ref{lemma:33} holds true. The Friesecke, James and M\"uller's rigidity estimate (Theorem 4.1 \cite{FJM02}) provides sequences of rotations that are piecewise constant on cubes of sizes $\varepsilon^2(h)$ with center in $\varepsilon^2(h)\ZZ^2$ but in this case the identification  of ``good cubes" and ``bad cubes" of size $\varepsilon^2(h)$ (Remark \ref{remarkcase1}) is not useful in this case because the contributions of the oscillations of the test functions over the cubes are not negligible. And for this case to characterize the three-scale limit $\tilde{V}$ for the flat domain $\Omega$
the condiction (\ref{conditiontoprovemach21}) becomes to
\begin{equation}\label{conditiontoprove2mach21}
\displaystyle\int_{\omega}\frac{(\tilde{R}^h(x')}{h}:(\nabla')^{\perp}\varphi^{\varepsilon}(x')\psi\left(x'\right)dx'=0,
\end{equation}
where the test function $\varphi^{\varepsilon}(x')=\varepsilon^2(h)\phi\left(\frac{x'}{\varepsilon(h)}\right)\varphi\left(\frac{x'}{\varepsilon^2(h)}\right)$ and the set
\begin{equation*}
\hat{\ZZ}^{\varepsilon}:=\left\{z\in\ZZ^2: Q(\varepsilon^2(h)z,\varepsilon^2(h))\cap\text{supp}\,\psi\neq \emptyset\right\},
\end{equation*}
plays  a crucial point to prove (\ref{conditiontoprove2mach21}).The fields $\mathcal{U}_{0, +\infty}$ arise, essentially, due the Lemma \ref{lemma:31}.\\

\textbf{Case 4:} $\gamma_1=0$ and $\gamma_2\in(0,+\infty)$.\\

We first apply Lemma \ref{lemma:33} with $\delta(h)=\varepsilon^2(h)$, and we can construct $\{\tilde{R}_{s}^h\}\subset L^{\infty}(\omega; SO(3))$ and $\{\tilde{R}^h\}\subset W^{1,2}(\omega;\RR^{3\times 3})$, satisfying (\ref{constructrotation}), and $\tilde{R}_{s}^h$ piecewise constant on every cube of the form
\begin{equation*}
Q(\varepsilon^2(h)z, \varepsilon^2(h)),\quad \text{with}\,z\in\ZZ^2.
\end{equation*}
Finally, arguing as in Case 3, we obtain the convergence desired.\\

\textbf{Case 5:} $\gamma_1=\gamma_2=0$.\\

Similarly to the previous case, we can apply Lemma \ref{lemma:33} with $\delta(h)=h$, and we can construct $\{\tilde{R}_{s}^h\}\subset L^{\infty}(\omega; SO(3))$ and $\{\tilde{R}^h\}\subset W^{1,2}(\omega;\RR^{3\times 3})$, satisfying (\ref{constructrotation}), and $\tilde{R}_{s}^h$ piecewise constant on every cube of the form
\begin{equation*}
Q(hz, h),\quad \text{with}\,z\in\ZZ^2.
\end{equation*}
We define,
\begin{equation*}
\ZZ^{\varepsilon^2}:=\left\{z\in\ZZ^2: Q(\varepsilon^2(h)z, \varepsilon^2(h))\times Q\cap\text{supp}\,\psi\neq \emptyset\right\}
\end{equation*}
and
\begin{equation*}
Q_{\varepsilon^2}:=\displaystyle\bigcup_{z\in \ZZ^{\varepsilon^2}}Q(\varepsilon^2(h) z, \varepsilon^2(h)).
\end{equation*}
We subdivide
\begin{equation*}
\mathcal{Q}_{h}:=\left\{Q(h\lambda,h): \lambda\in\ZZ^2\,\,\text{and}\,\, Q(h\lambda, h)\cap Q_{\varepsilon^2}\neq \emptyset\right\}
\end{equation*}
into two subsets:
\begin{itemize}
\item[(a)] ``good cubes of size $h$", i.e., those which are entirely contained in a cube of size $\varepsilon^2(h)$ belonging to $Q_{\varepsilon^2}$, and where $(R^h)'$ is hence constant,
\item[(b)] ``bad cubes of size $h$", i.e., those intersecting more than one element of $Q_{\varepsilon^2}$. 
\end{itemize}
In other words, we define the sets
\begin{equation*}
\ZZ^{h}_{g}:=\left\{\lambda\in\ZZ^2: \exists z\in\ZZ^{\varepsilon^2}\,\,\text{s.t.}\,\, Q(h\lambda, h)\subset Q(\varepsilon^2(h)z, \varepsilon^2(h))\right\}
\end{equation*}
and
\begin{equation*}
\ZZ^{h}_{b}:=\left\{\lambda\in\ZZ^2: Q(h\lambda, h)\cap Q_{\varepsilon^2}\neq \emptyset\quad\text{and}\quad\lambda\notin \ZZ^{h}_{g}\right\}
\end{equation*}
(where ``g" and ``b" stand for ``good" and ``bad", respectively). Finally, arguing as in the case 1, we obtain the characterization of $G^h$. The fields $\mathcal{U}_{0, 0}$ arise, essentially, due the Lemma \ref{lemma:31}.

\end{proof}

\begin{obs}\label{remarknewmay12}
The case $h\sim \varepsilon^2$ (i.e., $\gamma_1=0$ and $\gamma_2\in(0,+\infty)$) is critical because the term $\displaystyle\frac{1}{\gamma_2}\varphi\textbf{S}(x)$ appear in the three-scale limit of $G^h$. This phenomenon does not occur in the case of plates, where $\textbf{S}=0$.
\end{obs}

\vspace{0.5 cm}

\begin{lemma}\label{lemma:april30}
Let $(G^h)\subset L^2(S^1; \RR^3)$ be such that $G^h\overset{3}{\rightharpoonup} G$ in $L^2(S^1\times\mathcal{Y}\times\mathcal{Y};\RR^{3\times 3})$. Then,
\begin{equation}\label{eq:conver1}
h^{-1}\left(\sqrt{(Id+hG^h)^T(Id+hG^h)}-I\right)\overset{3}{\rightharpoonup}\text{sym}\, G\quad\text{in}\,L^2(S^1\times\mathcal{Y}\times\mathcal{Y};\RR^{3\times 3}).
\end{equation}
\end{lemma}

\begin{lemma}\label{lemma:3}
Let $(y^h)\subset W^{1,2}(S^1;\RR^3)$, define $E^h:S^1\to\RR^{3\times 3}$ by (\ref{eq:34}) and $E$ be such that $E^h\overset{3}{\rightharpoonup} E$. Then we have
\begin{align*}
\displaystyle\lim_{h\to 0}&\inf \int_{S}\int_{I}\mathscr{Q}\left(x+tn(x), r(x)/\varepsilon, r(x)/\varepsilon^2, E^h(x+tn(x))\right)\,dt\,d\text{vol}_{S}(x)\\
&\displaystyle\geq \int_{S}\int_{I}\int_{\mathcal{Y}}\int_{\mathcal{Y}}\mathscr{Q}\left(x+tn(x), y, z, E(x+tn(x))\right)\,dz\,dy\,dt\,d\text{vol}_{S}(x)\\
&\text{and}\\
\displaystyle\lim_{h\to 0}&\inf\frac{1}{h^2} \int_{S}\int_{I}W\left(x+tn(x), r(x)/\varepsilon, r(x)/\varepsilon^2, I+hE^h(x+tn(x))\right)\,dt\,d\text{vol}_{S}(x)\\
&\displaystyle\geq \int_{S}\int_{I}\int_{\mathcal{Y}}\int_{\mathcal{Y}}\mathscr{Q}\left(x+tn(x), y, z, E(x+tn(x),y,z)\right)\,dz\,dy\,dt\,d\text{vol}_{S}(x).
\end{align*}
\end{lemma}

\begin{proof}
    
 We refer to \cite{HORVEL15} and \cite{REF33} for a proof of Lemma \ref{lemma:3} in the case in which $\mathscr{Q}$ is independent of $z$. The proof in our setting is an adaptation.
\end{proof}
 \textbf{\underline {Proof of the Theorem \ref{teo:187}(i)} :} We first consider the case when $\gamma_1\in (0,+\infty)$ and $\gamma_2=+\infty$. Given the previous tools, we can provide a proof of the Theorem \ref{teo:187}(i). The proof of the lower bound follows standard arguments: truncation, Taylor expansion, and lower semicontinuity of integral functional (see Lemma \ref{lemma:3}) concerning the three-scale convergence of the sequence $E^h$, with three-scale limit $E=\text{sym}\,G$ (see Lemma \ref{lemma:april30}) where $G$ is obtained by Proposition \ref{prop:5}, i.e.
 \begin{align*}
&\displaystyle\liminf_{h\to 0}\frac{J^h(u^h)}{h^2}\\
&=\liminf_{h\to 0}\displaystyle\frac{1}{h^3}\int_{S^1}W(x,r(x)/\varepsilon, r(x)/\varepsilon^2, \nabla_{h}
y^h(x))\left|\det\nabla(\Theta^h)^{-1}(x)\right|dx\\
&\geq\liminf_{h\to 0}\displaystyle\frac{1}{h^2}\displaystyle\int_{S}\int_{I}W(x+tn(x),r(x)/\varepsilon, r(x)/\varepsilon^2,\nabla_{h}y^h(x+tn(x))\,dt\,d\,vol_{S}(x)\\
&= \liminf_{h\to 0}\displaystyle\frac{1}{h^2}\displaystyle\int_{S}\int_{I}W(x+tn(x),r(x)/\varepsilon, r(x)/\varepsilon^2,R^h\sqrt{(\nabla_{h}y^h(x+tn(x)))^T(\nabla_{h}y^h(x+tn(x)))}\,dt\,d\,vol_{S}(x)\\
&=\liminf_{h\to 0}\displaystyle\frac{1}{h^2}\displaystyle\int_{S}\int_{I}W(x+tn(x),r(x)/\varepsilon, r(x)/\varepsilon^2,\sqrt{(\nabla_{h}y^h(x+tn(x)))^T(\nabla_{h}y^h(x+tn(x)))}\,dt\,d\,vol_{S}(x)\\
&=\liminf_{h\to 0}\displaystyle\frac{1}{h^2}\displaystyle\int_{S}\int_{I}W(x+tn(x),r(x)/\varepsilon, r(x)/\varepsilon^2,Id+hE^h(x+tn(x)))\,dt\,d\,vol_{S}(x)\\
&\geq \int_{S}\int_{I}\int_{\mathcal{Y}}\int_{\mathcal{Y}}\mathscr{Q}\left(x+tn(x), y, z, E(x+tn(x),y,z)\right)\,dz\,dy\,dt\,d\text{vol}_{S}(x). 
\end{align*}
By applying Lemma \ref{lemma:3}, we get
\begin{align*}
&\int_{S}\int_{I}\int_{\mathcal{Y}}\int_{\mathcal{Y}}\mathscr{Q}\left(x+tn(x), y, z, E(x+tn(x),y,z)\right)\,dz\,dy\,dt\,d\text{vol}_{S}(x)\\
&= \displaystyle\int_{S}\int_{I}\int_{\mathcal{Y}}\int_{\mathcal{Y}}\mathscr{Q}\left(x+tn(x), y, z, \text{sym}\,G(x+tn(x),y,z)\right)\,dz\,dy\,dt\,d\text{vol}_{S}(x)\\
&=\displaystyle\int_{S}\int_{I}\int_{\mathcal{Y}}\int_{\mathcal{Y}}\mathscr{Q}\left(x+tn(x), y, z, \text{sym}\,B+t\textbf{S}^{r}_{u}+\mathcal{U}_{\gamma_1, +\infty}(\zeta, \eta, \rho)(x+tn(x),y,z)\right)\,dz\,dy\,dt\,d\text{vol}_{S}(x)\\
&\geq \inf_{U\in L_{\gamma_1,+\infty}(I\times \mathcal{Y}\times\mathcal{Y})}\displaystyle\int_{S}\int_{I}\int_{\mathcal{Y}}\int_{\mathcal{Y}}\mathscr{Q}\left(x+tn(x), y, z, \text{sym}\,B+t\textbf{S}^{r}_{u}+U(x+tn(x),y,z)\right)\,dz\,dy\,dt\,d\text{vol}_{S}(x)\\
&=\displaystyle\int_{S}\mathscr{Q}_{\gamma_1, +\infty}(x,\textbf{S}^{r}_{u}(x))\,d\text{vol}_{S}(x)
\end{align*}
For the other cases, $\gamma_1 = +\infty$ and $\gamma_2=+\infty$, $\gamma_1=0$ and $\gamma_2=+\infty$, $\gamma_1=0$ and $\gamma_2=(0,+\infty)$ and $\gamma_1=\gamma_2=0$, we use the same arugument with Proposition \ref{prop:5}.
\qed
\subsection{Proof of Theorem \ref{teo:187} ii). The upper bound}\label{ubound}
Let us introduce the recovery sequence. Recall Lemma 3.5 in \cite{HORVEL18}.
\begin{lemma}\label{lema:35}
Let $u\in W^{2,\infty}_{\text{iso}}(S)$ and define $\nu : S\to\mathbb{S}^2$ by
\begin{equation*}
\nu=\displaystyle\frac{\nabla_{\tau_1}u\times \nabla_{\tau_2}u}{|\nabla_{\tau_1}u\times \nabla_{\tau_2}u|}.
\end{equation*}
Let $w\in W^{2,\infty}(S,\RR^3)$ and define $\mu\in W^{1,\infty}(S,\RR^3)$ by
\begin{equation*}
\mu=\left(\nu \cdot \nabla_{\tau_1}w\right)\nabla_{\tau^1}u+\left(\nu \cdot \nabla_{\tau_2}w\right)\nabla_{\tau^2}u
\end{equation*}
and define the deformations $v^h: S^h\to\RR^3$ by
\begin{equation*}
v^h =u+t\nu+h(w+t\mu).
\end{equation*}

Define $R\in W^{1,\infty}(S, SO(3))$ by $R=\nabla uT_{S}+\nu\otimes n$. Then there exists $Y^h\in L^{\infty}(S^h, \RR^{3\times 3})$ with $||Y^h||_{L^{\infty}(S^h)}\leq Ch^2$ such that
\begin{equation*}
d v^h\odot R=I+h du \odot dw+t\textbf{S}^{r}_{u}+Y^h.
\end{equation*}
\end{lemma}

\begin{obs}
 The choice of our recovery sequences depends on the following two factors: 
\begin{itemize}
\item[(i)] It considers the inhomogeneity of material.
\item[(ii)] The energy density contains a spatial variable which makes it necessary to choose a nonzero displacement $w$ in Lemma \ref{lema:35} and whose existence is guaranteed by Proposition $2.15 $ in \cite{HORVEL18}.
\end{itemize}

Moreover, multilayered materials can be deduced as particular cases of Theorem \ref{teo:187}(cf. \cite{BUFDAVFON15} for the corresponding problem for plates).

\end{obs}

\textbf{\underline {Proof of the Theorem \ref{teo:187}(ii)} :}  By approximation, it is enough to prove the claim for $u\in W^{2,\infty}_{\text{iso}}(S)$ and, thanks also to Proposition 2.15 in \cite{HORVEL18}, for all $B$ of the form $B=du\odot dw$ with $w\in W^{2,\infty}(S,\RR^3)$.\\

\medskip We will use the same notation as in the statement of Lemma \ref{lema:35}; in particular, the definition of $v^h$ in terms of $w$ and $u$. Moreover, we set $\sigma^{\alpha}=\nabla_{\tau^{\alpha}}u$.

\vspace{0.5 cm}

\textbf{Case 1:} $\gamma_1\in (0,+\infty)$ and $\gamma_2=+\infty$. Let $\zeta\in C^1_{0}(S, \dot{C}^1(I\times\mathcal{Y},\RR^2))$, $\rho\in C^1_0(S, \dot{C}^1(I; \dot{C}^1(\mathcal{Y})))$ and $\eta\in C^1_0(S\times \mathcal{Y}, \dot{C}^1(I\times\mathcal{Y}))$ and define the rescaled deformations $y^h: S^1\to\RR^3$ by the following equation on $S^h$:
\begin{align*}
y^h(\Theta^h)&= v^h+h\varepsilon\zeta_{\alpha}\left(\pi, \frac{t}{h}, \frac{r}{\varepsilon}\right)\sigma^{\alpha}+h\varepsilon\rho\left(\pi, \frac{t}{h}, \frac{r}{\varepsilon}\right)\nu+h\varepsilon^2\eta_{\alpha}\left(\pi, \frac{t}{h}, \frac{r}{\varepsilon}, \frac{r}{\varepsilon^2}\right)\sigma^{\alpha}\\
&\quad+h\varepsilon^2\eta_{3}\left(\pi, \frac{t}{h}, \frac{r}{\varepsilon}, \frac{r}{\varepsilon^2}\right)\nu.
\end{align*}
Lemma \ref{lema:35} implies that on $S^1$
\begin{equation}\label{eq:2087}
\text{sym}\,(R^{T}\nabla_{h}y^h)=I+hB+th\textbf{S}^r_{u}+h\,\mathcal{U}_{\gamma_1, +\infty}(\zeta, \eta, \rho)\left(x, \frac{t}{h},\frac{r}{\varepsilon}, \frac{r}{\varepsilon^2}\right)+o(h),
\end{equation}
where $\lim_{h\to 0}||\frac{o(h)}{h}||=0$.

\medskip By frame invariance of $W$ and using (\ref{eq:linW}), we deduce from (\ref{eq:2087}) that
\begin{equation*}
\displaystyle\frac{1}{h^2}W\left(\cdot, \frac{r}{\varepsilon}, \frac{r}{\varepsilon^2}, \nabla_{h}y\right)\to \mathscr{Q}\left(\cdot, \frac{r}{\varepsilon}, \frac{r}{\varepsilon^2}, \textbf{S}^r_{u}+B+\mathcal{U}_{\gamma_1, +\infty}(\zeta, \eta, \rho)\left(\cdot, t,\frac{r}{\varepsilon}, \frac{r}{\varepsilon^2}\right)\right),
\end{equation*}
pointwise on $S^1$. From this we readily deduce
\begin{equation*}
\displaystyle\lim_{h\to 0}h^{-2}I^h (y^h)=\int_{S}\int_{I\times\mathcal{Y}\times\mathcal{Y}}\mathscr{Q}\left(\cdot, y, z, \textbf{S}^r_{u}+B+\mathcal{U}_{\gamma_1, +\infty}(\zeta, \eta, \rho)\left(\cdot, t, y, z\right)\right)\, dz\,dy\,dt\,d\,\text{vol}_{S}.
\end{equation*}

\textbf{Case 2:} $\gamma_1=\gamma_2=+\infty$. This is similar to the previous case. So, we only state the formula for the recovery sequence. For  $\zeta\in C^1_{0}(S, \dot{C}^1(I\times\mathcal{Y},\RR^2))$, $\rho\in C^1_0(S, \dot{C}^1(I; \dot{C}^1(\mathcal{Y})))$, $\eta\in C^1_0(S\times \mathcal{Y}, \dot{C}^1(I\times\mathcal{Y}))$ and $c\in C^1_0(S, C^1_0(I,\RR^3))$, we define $y^h: S^1\to\RR^3$ by the following equation on $S^h$:
\begin{align*}
y^h(\Theta^h)&= v^h+h\varepsilon\zeta_{\alpha}\left(\pi, \frac{t}{h}, \frac{r}{\varepsilon}\right)\sigma^{\alpha}+h\varepsilon\rho\left(\pi, \frac{t}{h}, \frac{r}{\varepsilon}\right)\nu+h\varepsilon^2\eta_{\alpha}\left(\pi, \frac{t}{h}, \frac{r}{\varepsilon}, \frac{r}{\varepsilon^2}\right)\sigma^{\alpha}\\
&\quad+h\varepsilon^2\eta_{3}\left(\pi, \frac{t}{h}, \frac{r}{\varepsilon}, \frac{r}{\varepsilon^2}\right)\nu+\displaystyle 2h^2\left(\int_{0}^{t/h}c_{\alpha}(x,s)\,ds\right)\sigma^{\alpha}\\
&\quad+h^2\left(\int_{0}^{t/h}c_3(x,s)\,ds\right)\nu.
\end{align*}
By Lemma \ref{lema:35}, the frame invariance of $W$ and similar computations of case 1, we obtain:
\begin{equation*}
\displaystyle\lim_{h\to 0}h^{-2}I^h (y^h)=\int_{S}\int_{I\times\mathcal{Y}\times\mathcal{Y}}\mathscr{Q}\left(\cdot, y, z, \textbf{S}^r_{u}+B+\mathcal{U}_{+\infty, +\infty}(\zeta,\eta, \rho, c)\left(\cdot, t, y, z\right)\right)\, dz\,dy\,dt\,d\,\text{vol}_{S}.
\end{equation*}

\textbf{Case 3:} $\gamma_1=0$ and $\gamma_2=+\infty$. For  $\zeta\in C^1_{0}(S, \dot{C}^1(\mathcal{Y},\RR^2))$, $\varphi\in C^2_0(S, \dot{C}^2(\mathcal{Y}))$, $\eta\in C^1_0(S\times \mathcal{Y}, \dot{C}^1(I\times\mathcal{Y}))$ and $\mu \in C^1_0(S, C^1_0(I\times\mathcal{Y},\RR^3))$, we define $y^h: S^1\to\RR^3$ by the following equation on $S^h$:
\begin{align*}
y^h(\Theta^h)&= v^h+h\varepsilon\zeta_{\alpha}\left(\pi, \frac{r}{\varepsilon}\right)\sigma^{\alpha}+\varepsilon^2\varphi\left(\pi, \frac{r}{\varepsilon}\right)\nu-t\varepsilon\partial_{y_{\alpha}}\varphi\left(\pi, \frac{r}{\varepsilon}\right)\sigma^{\alpha}-t\varepsilon^2\partial_{\alpha}\varphi\left(\pi, \frac{r}{\varepsilon}\right)\sigma^{\alpha}\\
&\quad\displaystyle+h\varepsilon^2\eta_{\alpha}\left(\pi, \frac{t}{h}, \frac{r}{\varepsilon}, \frac{r}{\varepsilon^2}\right)\sigma^{\alpha}+h\varepsilon^2\eta_3\left(\pi, \frac{t}{h}, \frac{r}{\varepsilon}, \frac{r}{\varepsilon^2}\right)\nu\\
&\quad+\displaystyle 2h^2\left(\int_{0}^{t/h}\mu_{\alpha}\left(\pi,s, \frac{r}{\varepsilon}\right)\,ds\right)\sigma^{\alpha}+h^2\left(\int_{0}^{t/h}\mu_3\left(\pi, s, \frac{r}{\varepsilon}\right)\,ds\right)\nu.
\end{align*}
By Lemma \ref{lema:35}, the frame invariance of $W$ and similar computations of case 1, we obtain:
\begin{equation*}
\displaystyle\lim_{h\to 0}h^{-2}I^h (y^h)=\int_{S}\int_{I\times\mathcal{Y}\times\mathcal{Y}}\mathscr{Q}\left(\cdot, y, z, \textbf{S}^r_{u}+B+\mathcal{U}_{0, +\infty}(\zeta,\eta, \varphi, \mu)\left(\cdot, t, y, z\right)\right)\, dz\,dy\,dt\,d\,\text{vol}_{S}.
\end{equation*}

\textbf{Case 4:}  $\gamma_1=0$ and $\gamma_2\in(0,+\infty)$. For  $\zeta\in C^1_{0}(S, \dot{C}^1(\mathcal{Y},\RR^2))$, $\varphi\in C^2_0(S, \dot{C}^2(\mathcal{Y}))$, $\eta\in C^1_0(S\times \mathcal{Y}, \dot{C}^1(I\times\mathcal{Y}))$ and $\mu \in C^1_0(S, C^1_0(I\times\mathcal{Y},\RR^3))$, we define $y^h: S^1\to\RR^3$ by the following equation on $S^h$:
\begin{align*}
y^h(\Theta^h)&= v^h+h\varepsilon\zeta_{\alpha}\left(\pi, \frac{r}{\varepsilon}\right)\sigma^{\alpha}+\varepsilon^2\varphi\left(\pi, \frac{r}{\varepsilon}\right)\nu-t\varepsilon\partial_{y_{\alpha}}\varphi\left(\pi, \frac{r}{\varepsilon}\right)\sigma^{\alpha}-t\varepsilon^2\partial_{\alpha}\varphi\left(\pi, \frac{r}{\varepsilon}\right)\sigma^{\alpha}\\
&\quad\displaystyle+h\varepsilon^2\eta_{\alpha}\left(\pi, \frac{t}{h}, \frac{r}{\varepsilon}, \frac{r}{\varepsilon^2}\right)\sigma^{\alpha}+h\varepsilon^2\eta_3\left(\pi, \frac{t}{h}, \frac{r}{\varepsilon}, \frac{r}{\varepsilon^2}\right)\nu\\
&\quad+\displaystyle 2h^2\left(\int_{0}^{t/h}\mu_{\alpha}\left(\pi,s, \frac{r}{\varepsilon}\right)\,ds\right)\sigma^{\alpha}+h^2\left(\int_{0}^{t/h}\mu_3\left(\pi, s, \frac{r}{\varepsilon}\right)\,ds\right)\nu.
\end{align*}
Finally, by Lemma \ref{lema:35}, the frame invariance of $W$ and similar computations of case 1, we obtain:
\begin{equation*}
\displaystyle\lim_{h\to 0}h^{-2}I^h (y^h)=\int_{S}\int_{I\times\mathcal{Y}\times\mathcal{Y}}\mathscr{Q}\left(\cdot, y, z, \textbf{S}^r_{u}+B+\mathcal{U}_{0, \gamma_2}(\zeta,\eta, \varphi, \mu)\left(\cdot, t, y, z\right)\right)\, dz\,dy\,dt\,d\,\text{vol}_{S}.
\end{equation*}
\lqqd

\section{$\Gamma-$limit for convex shells}\label{sec:6}

In this section, we shall identify the $\Gamma-$limit for convex shells in the remaining case $\gamma_1=\gamma_2=0$, \textit{i.e.} $h\ll\varepsilon^2$. We wish to illustrate the stronger influence of the geometry in this case. For obtaining the limit model we shall closely follow the arguments used in \cite{HORVEL15} as follows: We work under the assumption that $S$ is uniformly convex, $i.e.$, there exists $C>0$ such that
\begin{equation}\label{eq:6001}
\textbf{S}(x)\tau\cdot\tau\geq C|\tau|^2_{T_xS},\quad\forall x\in S, \tau\in T_xS.
\end{equation}
For $x\in S$ we define a relaxation operator with the values in $L^2(I\times\mathcal{Y}\times\mathcal{Y};\RR^{3\times3}_{\text{sym}})$ as follows: Set $D(\mathcal{U}_{0,0})=\dot{L}^2(\mathcal{Y};\RR^{2\times 2}_{\text{sym}})\times L^2(I\times\mathcal{Y}; \dot{W}^{1,2}(\mathcal{Y}; \RR^3))\times L^2(I\times\mathcal{Y};\RR^{3})$ and for all $(\dot{B}, \eta, \mu)\in L^2(S, D(\mathcal{U}_{0, 0}))$ define
\begin{align*}
\mathcal{U}_{0, 0}(\dot{B},\eta,\mu)&=\displaystyle\sum_{i,j=1}^{3}\begin{pmatrix}
  \dot{B}  & \begin{matrix}\,\,\mu_1 \\
  \,\,\mu_2 
  \end{matrix} \\
 (\mu_1,\mu_2) & \mu_3
\end{pmatrix}_{ij}\tau^{i}\otimes\tau^j+(\text{sym}\,\nabla_{z}\eta|0)_{ij}\tau^{i}\otimes\tau^j.
\end{align*}
As usual, we introduce the vector bundle $L_{0, 0}(I\times\mathcal{Y}\times\mathcal{Y})$ of relaxation fields to be the range of $\mathcal{U}_{0,0}$ similarly to the bundles $L_{0,+\infty}(I\times\mathcal{Y}\times\mathcal{Y})$ introduced earlier. As in the previous cases, each fiber of $L_{0, 0}(I\times\mathcal{Y}\times\mathcal{Y})$ is a closed subspace of $L^2(I\times\mathcal{Y}\times\mathcal{Y};\RR^{3\times3}_{\text{sym}})$. We also define the functional $I_{0,0}: W^{1,2}(S;\RR^3)\to\RR$ by setting
\begin{equation}
I_{0,0}(u)=\begin{cases}
\int_{S}\mathscr{Q}_{0,0}(x, \textbf{S}^{r}_{u}(x))\,d\,\text{vol}_{S}(x) &\text{if}\, u\in W^{2,2}_{\text{iso}}(S), \\
+\infty &\text{otherwise},
\end{cases}
\end{equation}
with the quadratic form $\mathscr{Q}_{0,0}(x, \cdot):T^{*}S\odot T^{*}S\to\RR$ given by
\begin{equation}\label{eq:6003}
\mathscr{Q}_{0,0}(x, q)=\inf \int_{I}\int_{\mathcal{Y}}\int_{\mathcal{Y}}\mathscr{Q}\left(x+tn(x), y, z, p+tq+U(t,y,z)\right)\,dy\,dz\,dt.
\end{equation}
Here the infimum is taken over all $U\in L^{(x)}_{0,0}(I\times\mathcal{Y}\times\mathcal{Y})$ and all $p\in T^{*}_{x}S\otimes T^{*}_{x}S$.
\vspace{0.25 cm}

We introduce the space
\begin{align*}
&FL(S; \dot{C}^{\infty}(\mathcal{Y}))\\
&=\left\{(x,y)\mapsto \displaystyle\sum_{k\in \ZZ^2,\,|k|\leq n,\, k\neq 0}c^k(x)\, e^{2\pi i k\cdot y}: n\in\NN\,\,\text{and}\,\, c^k\in C^{1}_{0}(S; \CC)\,\,\text{with}\,\, \overline{c}^k=c^{-k}\right\}.
\end{align*}
By Fourier transform it can be easily seen that $FL(S; \dot{C}^{\infty}(\mathcal{Y}))$ is dense in $L^2(S; \dot{H}^m(\mathcal{Y}))$, for any $m\in\NN_{0}$.

\vspace{0.25 cm}

Let us recall the following result that appears in Proposition 4.2, \cite{HORVEL15}.

\begin{prop}\label{prop:3horvel}
Let $(w^h)$ be a bounded sequence in $H^2(S;\RR^3)$ such that $\frac{1}{h}q_{w^h}$ is bounded in $L^2(S; T^{*}S\otimes T^{*}S)$. Then there exist $w_0\in H^2(S)$, $w_1\in L^2(S; \dot{H}^2(\mathcal{Y};\RR^3))$ and $B\in L^2(S, \dot{L}^2(\mathcal{Y}; T^{*}S\otimes T^{*}S))$ such that, after passing to a subsequence, $q_{w^h}/h\overset{2}{\rightharpoonup} B$ and $\text{Hess}\,w^h\overset{2}{\rightharpoonup}\text{Hess}\,w_0+\text{Hess}_{\mathcal{Y}}\,w_1$. Set $B_{w}=\int_{Y}B(\cdot, y)\,dy$. Then the following condition are true:
\begin{itemize}
\item[(i)] If $h\gg \varepsilon^2$ then there exists a unique $v\in L^2(S; \dot{H}^{1}(\mathcal{Y}; \RR^2))$ such that
\begin{equation*}
B=B_{w}+\text{Def}_{\mathcal{Y}}\,v.
\end{equation*}

\item[(ii)] If $h\sim \varepsilon^2$ and if we set $\frac{1}{\gamma_2}=\lim_{h\to 0}\frac{\varepsilon^2}{h}$,  then there exists a unique $v\in L^2(S; \dot{H}^{1}(\mathcal{Y}; \RR^2))$ such that
\begin{equation*}
B=B_{w}+\text{Def}_{\mathcal{Y}}\,v+\displaystyle\frac{1}{\gamma_2}(w_1\cdot n)\textbf{S}.
\end{equation*}

\item[(iii)] If $h\ll \varepsilon^2$, then there exists a unique $v\in L^2(S; \dot{H}^{1}(\mathcal{Y}; \RR^2))$ such that
\begin{equation*}
\text{Def}_{\mathcal{Y}}\,v+\displaystyle(w_1\cdot n)\textbf{S}=0.
\end{equation*}
\end{itemize}
\end{prop}

\begin{lemma}\label{lemma:46}(see Lemma 6.1 in \cite{HORVEL15})
Assume (\ref{eq:6001}) is satisfied and let $\dot{B}\in L^2(S; \dot{L}^2(\mathcal{Y}; T^{*}S\otimes T^{*}S))$. Then there exists unique $w\in L^2(S; \dot{H}^1(\mathcal{Y};\RR^2))$ and $\varphi\in L^2(S; \dot{L}^2(\mathcal{Y}))$ such that
\begin{equation}
\text{Def}_{\mathcal{Y}}w+\varphi\textbf{S}=\dot{B}.
\end{equation}
Moreover, if $\dot{B}_{ij}\in FL(S; \dot{C}^{\infty}(\mathcal{Y}))$ for every $i,j=1, 2$ then $w_i\in FL(S; \dot{H}^1(\mathcal{Y}))$, for $i=1, 2$ and $\varphi\in FL(S; \dot{H}^1(\mathcal{Y}))$.
\end{lemma}

\begin{teo}\label{theo61}
Under the hypotheses and with the notation of Theorem \ref{teo:187} and assuming, in addition, that $S$ is uniformly convex and that $h\ll \varepsilon^2$, moreover, the following are true:
\begin{itemize}
\item[$\bullet$] We have
\begin{equation*}
\lim_{h\to 0}\inf h^{-2} J^h(u^h)\geq I_{0,0}(u).
\end{equation*}
\item[$\bullet $] If, in addition, $S$ is simply connected, then for every $u\in \widetilde{W}^{2, 2}_{\text{iso}}(S)$ there exists $(u^h)\subset W^{1,2}(S^h;\RR^3)$ satisfying (\ref{eq:33}) and such that $y^h\to u$, strongly in $W^{1,2}(S^1)$. Moreover,
\begin{equation*}
\lim_{h\to 0}h^{-2}J^h(u^h)=I_{0,0}(u).
\end{equation*}
\end{itemize}
\end{teo}

\begin{proof}
We only sketch the proof. As in Proposition \ref{prop:5} there exist $B\in L^2(S, T^{*}S\odot T^{*}S)$ and $(\zeta, \eta, \varphi, \mu)\in L^2(S, D(\mathcal{U}_{0, +\infty}))$ such that (\ref{eq:35}) is satisfied. Using Proposition \ref{prop:3horvel} (iii) as well Lemma \ref{lemma:46}, we conclude that $\varphi=0$. Thus by Proposition \ref{prop:5} there exists $(\dot{B}, \eta,\mu)\in L^2(S, D(\mathcal{U}_{0,0}))$ , where $\dot{B}=\text{Def}_{\mathcal{Y}}\zeta$ and $B\in L^2(S, T^{*}S\odot T^{*}S)$ such that the maps $G^h$ defined as in (\ref{eq:34}) converge weakly three-scale to
\begin{equation*}
 G=B+t\textbf{S}^{r}_{u}+\mathcal{U}_{0, 0}(\dot{B}, \eta, \mu).
\end{equation*}
Hence the lower bound part follows readily from the Lemma \ref{lemma:3} and definition of the functional $I_{0,0}$.

To prove the upper bound part we consider $\dot{B}$ with $(\dot{B})_{ij}\in FL(S;\dot{C}^{\infty}(\mathcal{Y}))$ for $i,j=1,2$. From Lemma \ref{lemma:46} there exists $\zeta\in(FL(S;\dot{C}(\mathcal{Y})))^2$ and $\varphi\in FL(S;\dot{C}^{\infty}(\mathcal{Y}))$ solving the system $\text{Def}_{\mathcal{Y}}\zeta+\varphi\textbf{S}=\dot{B}$. We choose $\eta\in C^1_0(S\times \mathcal{Y}, \dot{C}^1(I\times\mathcal{Y}))$ and $\mu \in C^1_0(S, C^1_0(I\times\mathcal{Y},\RR^3))$ and we define $y^h: S^1\to\RR^3$ by the following equation on $S^h$:
\begin{align*}
y^h(\Theta^h)&= v^h+h\varepsilon\zeta_{\alpha}\left(\pi, \frac{r}{\varepsilon}\right)\sigma^{\alpha}+h\varphi\left(\pi, \frac{r}{\varepsilon}\right)\nu\\
&\quad\displaystyle+h\varepsilon^2\eta_{\alpha}\left(\pi, \frac{t}{h}, \frac{r}{\varepsilon}, \frac{r}{\varepsilon^2}\right)\sigma^{\alpha}+h\varepsilon^2\eta_3\left(\pi, \frac{t}{h}, \frac{r}{\varepsilon}, \frac{r}{\varepsilon^2}\right)\nu\\
&\quad+\displaystyle 2h^2\left(\int_{0}^{t/h}\mu_{\alpha}\left(\pi,s, \frac{r}{\varepsilon}\right)\,ds\right)\sigma^{\alpha}+h^2\left(\int_{0}^{t/h}\mu_3\left(\pi, s, \frac{r}{\varepsilon}\right)\,ds\right)\nu.
\end{align*}

Now, we can argue as in the proof of upper bound in Section \ref{ubound} to conclude the desired result.
\end{proof}

\begin{appendices}
\section{Appendix}
In this section, for the sake of convenience of interested readers, we compile some results without proof and with proper references which play important roles in the proof of our Lemma  \ref{lemma:33} and Proposition \ref{prop:5} for the construction of the rotations $R^h$ that are piecewise constant on cubes of the form $Q(\varepsilon(h)z, \varepsilon(h))$. To be more precise, among the resuls below, Lemma \ref{lemma:a1} and Lemma \ref{lemma:a2} are helpful for any cases where as Lemma \ref{Lemma:a3} is specifically for $\gamma_1=\gamma_2=+\infty$. 
\begin{lemma}\label{lemma:a1}(See Lemma 3.1 in \cite{BUFDAVFON15})
Let $\gamma\in (0,1]$ and let $h, \delta >0$ be such that
\begin{equation*}
\gamma_0\leq \displaystyle\frac{h}{\delta}\leq \displaystyle\frac{1}{\gamma_0}.
\end{equation*}
There exists a constant $C$, depending only on $\omega$ and $\gamma_0$, such that for every $u\in W^{1,2}(\omega;\RR^3)$ there exists a map $R:\omega\to SO(3)$ piecewise constant on each cube $x+\delta Q$, with $x\in \delta \ZZ^2$, and there exists $\tilde{R}\in W^{1,2}(\omega;\RR^{3\times 3})$ such that
\begin{align}
&\nonumber||\nabla_h u-R||^2_{L^2(\Omega;\RR^3)}+||R-\tilde{R}||^2_{L^2(\omega;\RR^3)}+h^2||\nabla'\tilde{R}||^2_{L^2(\omega;\RR^3\times\RR^3)}\\
&\label{rotationRandRtilde}\quad\quad\leq C||\text{dist}\,(\nabla_h u;SO(3))||_{L^2(\Omega)}.
\end{align}
Moreover, for every $\xi\in\RR^2$ satisfying
\begin{equation*}
|\xi|_{\infty}:=\max \left\{|\xi\cdot e_1|, |\xi\cdot e_2|\right\}<h,
\end{equation*}
and for every $\omega'\subset \omega$, with $\text{dist}\,(\omega',\partial \omega)>Ch$, there holds
\begin{equation*}
||R(x')-R(x'+\xi)||_{L^2(\omega';\RR^2)}\leq C||\text{dist}\,(\nabla_h u;SO(3))||^2_{L^2(\omega)}.
\end{equation*}
\end{lemma}
\begin{lemma}\label{lemma:a2}(See Lemma 3.3 in \cite{HORVEL18})
Define
\begin{equation*}
\delta=
\begin{cases}
\varepsilon, & \text{if}\,\, \gamma_1\in(0,+\infty), \gamma_2=+\infty,\\
\displaystyle\left\lceil\frac{h}{\varepsilon}\right\rceil\varepsilon, &\text{if}\,\,\gamma_1=\gamma_2=+\infty,\\
h, & \text{if}\,\,\gamma_1=0, \gamma_2=+\infty.
\end{cases}
\end{equation*}
Then there exist constants $C, c>0$ such that the following is true: if $h\leq c$ and $u\in W^{1,2}(S^h;\RR^3)$, then there exists a map $\tilde{R}:\omega\to SO(3)$ which is constant on each cube $x+\delta Y$ with $x\in\delta \ZZ$ and there exist $\tilde{R}_{s}\in W^{1,2}(\omega;\RR^3)$ such that for each $a\in\RR^2$ with $|a_1|\leq \delta$ and $|a_2|\leq \delta$ and for each $\tilde{\omega}\subset\omega$ with $\text{dist}\,(\tilde{\omega}, \partial\omega)>c\delta$ we have:
\begin{align*}
&||(\nabla_{h}y)(\Xi)-\tilde{R}||^2_{L^2(\tilde{\omega}\times I)}+||\tilde{R}-\tilde{R}_{s}||^2_{L^2(\tilde{\omega})}+h^2||\tilde{R}-\tilde{R}_{s}||^2_{L^{\infty}(\tilde{\omega})}\\
&+h^2||(\partial_{1}\tilde{R}_{s}, \partial_{2}\tilde{R}_{s})||^2_{L^2(\tilde{\omega})}+||\tilde{R}(\cdot+a)-\tilde{R}||^2_{L^2(\tilde{\omega})}\\
&\leq C\displaystyle\int_{\Omega}\text{dist}^2(\nabla_{h}y(\Xi), SO(3)).
\end{align*}
\end{lemma}

\vspace{0.5 cm}

The following results will be useful for proving the Proposition \ref{prop:5}: For every $z\in\ZZ^2$ there exists $z'\in\ZZ^2$ such that
\begin{equation*}
Q(\varepsilon(h)z, \varepsilon(h))\subset Q(\delta(h),z'\delta(h))
\end{equation*}
or equivalently, with $m=\displaystyle\frac{\delta(h)}{\varepsilon(h)}\in\NN$,
\begin{equation}\label{eq:a1}
 \displaystyle\left(z-\frac{1}{2},z+\frac{1}{2}\right)\subset m\left(z'-\frac{1}{2}, z'+\frac{1}{2}\right).   
\end{equation}

\begin{lemma}\label{Lemma:a3}(See Lemma A.1. in \cite{BUFDAVFON15})Let $a\in\NN_0$. Then for every $z\in\ZZ$ there exists $z'\in\ZZ$ such that (\ref{eq:a1}) holds with $m=2a+1$.
\end{lemma}
Now we mention the observation which is used in Lemma \ref{lemma:33} for the case $\gamma_1\in(0,+\infty)$ and $\gamma_2=+\infty$.
\begin{obs}\label{remarkcase1}
If $\psi\in C^{\infty}_{c}(\omega; C^{\infty}_{\text{per}}(Y))$ and $h\to 0$, we can assume, without loss of generality, that for $h$ small enough
\begin{equation*}
\text{dist}(\text{supp}\,\psi; \partial\omega\times Y)>\left(1+\displaystyle\frac{3}{\gamma_1}\right)h.
\end{equation*}
We define
\begin{equation*}
\ZZ^{\varepsilon}:=\left\{z\in \ZZ^2: Q\left(\varepsilon(h) z, \varepsilon (h)\right)\times Y\cap \text{supp}\,\psi \neq\emptyset\right\}
\end{equation*}
and
\begin{equation*}
Q_{\varepsilon}:=\displaystyle\bigcup_{z\in \ZZ^{\varepsilon}}Q(\varepsilon(h)z, \varepsilon(h)).
\end{equation*}
If $\gamma_1\in (0,+\infty)$, for $h$ small enough we have $\sqrt{2}\varepsilon (h)<\frac{2h}{\gamma_1}$, so that
\begin{equation*}
\text{dist}\,(Q_{\varepsilon};\partial \omega)\geq \left(1+\displaystyle\frac{3}{\gamma_1}\right)h-\sqrt{2}\varepsilon (h)\geq \left(1+\displaystyle\frac{1}{\gamma_1}\right)h.
\end{equation*}
We subdivide
\begin{equation*}
\mathcal{Q}_{\varepsilon^2}:=\left\{Q(\varepsilon^2(h)\lambda,\varepsilon^2(h)): \lambda\in\ZZ^2\,\,\text{and}\,\, Q(\varepsilon^2(h)\lambda, \varepsilon^2(h))\cap Q_{\varepsilon}\neq \emptyset\right\}
\end{equation*}
into two subsets:
\begin{itemize}
\item[(a)] ``good cubes of size $\varepsilon^2(h)$", i.e., those which are entirely contained in a cube of size $\varepsilon(h)$ belonging to $Q_{\varepsilon}$, and where $(R^h)'$ is hence constant,
\item[(b)] ``bad cubes of size $\varepsilon^2(h)$", i.e., those intersecting more than one element of $Q_{\varepsilon}$. 
\end{itemize}
We observe that, if $\gamma_2=+\infty$,
\begin{equation}\label{eq:Apendix1}
\text{dist}\,(\mathcal{Q}_{\varepsilon^2}; \partial\omega)\geq \text{dist}\,(Q_{\varepsilon};\partial\omega)-\sqrt{2}\varepsilon^2(h)>h
\end{equation}
for $h$ small enough, and
\begin{equation}\label{eq:Apendix2}
\#\ZZ^{\varepsilon}\leq \displaystyle C\frac{|\omega|}{\varepsilon^2(h)}.
\end{equation}
Moreover, if $z\in\ZZ^{\varepsilon}$, $\lambda\in\ZZ^2$, and
\begin{equation*}
\varepsilon^2(h)\lambda\in Q(\varepsilon(h) z, \varepsilon(h)-\varepsilon^2(h)),
\end{equation*}
then $Q(\varepsilon^2(h)\lambda, \varepsilon^2(h))$ is a ``good cube", therefore that the boundary layer of $Q(\varepsilon(h)z, \varepsilon(h))$, that could possibly intersect ``bad cubes" has measuring given by
\begin{align*}
&|Q(\varepsilon(h)z, \varepsilon(h))|-|Q(\varepsilon(h)z, \varepsilon(h)-\varepsilon^2(h))|\\
&\quad=\varepsilon^2(h)-\left(\varepsilon(h)-\varepsilon^2(h)\right)^2=2\varepsilon^3(h)-\varepsilon^4(h).
\end{align*}
By (\ref{eq:Apendix2}) we conclude that the sum of all areas of ``bad cubes" intersecting $Q_{\varepsilon}$ is bounded from above by
\begin{equation}\label{eq:Apendix3}
\displaystyle C \frac{|\omega|}{\varepsilon^2(h)}(2\varepsilon^3(h)-\varepsilon^4(h))\leq C\varepsilon(h).
\end{equation}
We define the sets
\begin{equation*}
\ZZ^{\varepsilon}_{g}:=\left\{\lambda\in\ZZ^2: \exists z\in\ZZ^{\varepsilon}\,\,\text{s.t.}\,\, Q(\varepsilon^2(h)\lambda, \varepsilon^2(h))\subset Q(\varepsilon(h)z, \varepsilon(h))\right\}
\end{equation*}
and
\begin{equation*}
\ZZ^{\varepsilon}_{b}:=\left\{\lambda\in\ZZ^2: Q(\varepsilon^2(h)\lambda, \varepsilon^2(h))\cap Q_{\varepsilon}\neq \emptyset\quad\text{and}\quad\lambda\notin \ZZ^{\varepsilon}_{g}\right\}
\end{equation*}
(where ``g" and ``b" stand for ``good" and ``bad", respectively).
\end{obs}
Here we present the fact which is used in Proposition \ref{prop:5} for the case $\gamma_1=0$,$\gamma_2=+\infty$.
\begin{obs}\label{remarkA2}
By Lemma \ref{Lemma:a3} it follows that, setting $p:=\displaystyle\frac{\delta(h)}{\varepsilon^2(h)}$ and provided $p$ is odd, for every $z\in\ZZ^2$ there exists $z'\in\ZZ^2$ such that
\begin{equation*}
 Q(\varepsilon^2(h)z, \varepsilon^2(h))\subset Q(\delta(h)z',\delta(h)).   
\end{equation*}
\end{obs}

\begin{obs}\label{remarkA3}
We point out that if $m$ is even there many be $z\in\ZZ$ such that (\ref{eq:a1}) fails to be true for $z'\in\ZZ$, i.e.
\begin{equation*}
 \displaystyle\left(z-\frac{1}{2},z+\frac{1}{2}\right)\nsubseteq m\left(z'-\frac{1}{2}, z'+\frac{1}{2}\right).   
\end{equation*}
Indeed, if $m$ is even, then $z=\displaystyle\frac{3}{2}m\in\NN$ and
\begin{equation*}
 \begin{cases}
    \displaystyle z-\frac{1}{2}\geq (2a+1)z'-\frac{(2a+1)}{2}, & \\
    & \\
    \displaystyle z+\frac{1}{2}\leq (2a+1)z'+\frac{(2a+1)}{2}, &
 \end{cases}   
\end{equation*}
which in turn is equivalent to
\begin{equation*}
z'\in\left[1+\displaystyle\frac{1}{2m}, 2-\frac{1}{2m}\right]. \end{equation*}
This last condition leads to contradiction as
\begin{equation*}
 \left[1+\displaystyle\frac{1}{2m}, 2-\frac{1}{2m}\right]\cap\ZZ=\emptyset,\,\,\text{for every}\,\,m\in\NN.    
\end{equation*}
\end{obs}

\vspace{0.5 cm}

\end{appendices}

\section*{Acknowledgment}
This work was inspired by the great contributions of Prof. Igor Vel\v ci\'c to the theory of elasticity for composite materials for thin structures. P. Hern\'andez-Llanos was supported by Agencia Nacional de Investigaci\'on y Desarro\-llo de Chile (ANID) through FONDECYT Postdoctorado 2023 Grant No. 3230202 and supported also by Croatian Science Foundation under Grant Agreement No. IP-2018-01-8904 (Homdirestroptcm). The warm hospitality at the Faculty of Electrical Engineering and Computing (FER) of the University of Zagreb are gratefully acknowledged by P. Hern\'andez-Llanos. The central part of this paper was done while the third author was affiliated as a researcher at FER.\\
The authors T. Durante and  L. Faella  are members of GNAMPA (INDAM) and have been partially supported by GNAMPA (INDAM) under the project GNAMPA 2024-CUP E53C23001670001 (Italy).

\end{document}